\documentclass[11pt, a4paper, reqno]{amsart}
\usepackage{amssymb,amsmath}
\usepackage[colorlinks=true,linkcolor=blue,citecolor=blue]{hyperref}
\usepackage[margin=1in]{geometry}
\usepackage{mlmodern}
\usepackage[T1]{fontenc}
\usepackage{mathrsfs}
\usepackage{comment}
\usepackage{xcolor}

\newcommand{\normm}[1]{|\!|\!|#1|\!|\!|}

\numberwithin{equation}{section}
\newtheorem{theorem}{Theorem}[section]
\newtheorem{lemma}[theorem]{Lemma}
\newtheorem{proposition}[theorem]{Proposition}

\newtheorem{definition}[theorem]{Definition}
\newtheorem{remark}[theorem]{Remark}
\newcommand{\wlim}{\operatorname{w-lim}}
\newcommand{\R}{\mathbb R}

\DeclareMathOperator{\diam}{diam}
\newcommand{\bra}[1]{\left(#1\right)}

\begin{document}
\title[Non-autonomous KdV-Burgers equations]
{Large time dynamics for non-autonomous weakly damped KdV-Burgers equations on the real line}

\author [C.T. Anh, N. Ayhan, B.Q. Tang, I. Witt] {Cung The Anh, Nes\.ibe Ayhan, Bao Quoc Tang$^{\natural}$ and
  Ingo Witt}

\address{Cung The Anh \hfill\break
Faculty of Mathematics and Informatics, Hanoi National University of Education \hfill\break
136 Xuan Thuy, Cau Giay, Hanoi, Vietnam}
\email{anhctmath@hnue.edu.vn} 

\address{Nes\.ibe Ayhan \hfill\break 
Department of Mathematics and Scientific Computing, University of Graz\hfill\break Heinrichstrasse 36, 8010 Graz, Austria
}
\email{nesibe.ayhan@uni-graz.at}

\address{Bao Quoc Tang \hfill\break 
Department of Mathematics and Scientific Computing, University of Graz\hfill\break Heinrichstrasse 36, 8010 Graz, Austria
}
\email{quoc.tang@uni-graz.at}

\address{Ingo Witt \hfill\break
Mathematical Institute, University of G\"{o}ttingen, 
Bunsenstra{\ss}e 3-5, D-37073 G\"{o}ttingen, Germany}
\email{iwitt@uni-math.gwdg.de} 

\thanks{$^{\natural}$ Corresponding author: quoc.tang@uni-graz.at}

\begin{abstract}
	The large time dynamics of the weakly damped non-autonomous Korteweg-de Vries-Burgers equation
\begin{equation*}
u_t + u_{xxx}+ uu_x + \lambda u - \mu  u_{xx} = g(x,t)
\end{equation*}
 in the phase space $L^2(\mathbb R)$ is investigated. Due to the low regularity of the initial data, the framework of Bourgain spaces is employed. We first obtain the existence of a unique pullback attractor $\{\mathcal A_\mu(t)\}_{t\in\mathbb R}$ in $L^2(\mathbb R)$ for each $\mu \ge 0$. Then we show an optimal asymptotic smoothing effect of the equation, namely the pullback attractor is compact in $H^3(\mathbb R)$. These results are proved by using an energy equation method, splitting solutions into parts with low and high frequencies and exploiting a dispersive regularization. Furthermore, the attractors are shown to be upper semicontinuous as $\mu \to 0$ and have finite fractal dimensions in $H^1(\mathbb R)$.
\end{abstract}

\maketitle

\tableofcontents

\section{Introduction}

This paper is devoted to the following non-autonomous weakly damped Korteweg-de Vries-Burgers
(KdV-Burgers for short) equation: for $\tau \in \mathbb R$,
\begin{equation}\label{e1}
\begin{cases} 
\enspace u_t + u_{xxx}+uu_{x}+\lambda u - \mu  u_{xx} = g(t,x), & x\in \mathbb R, \quad t>\tau,\\
\enspace u(\tau,x) = u_{\tau}(x), &x\in \mathbb R,
\end{cases}
\end{equation}
where $\lambda>0$ and $\mu \in (0,1]$, the initial data
$u_\tau$ belongs to $L^2(\mathbb R)$, and the time-dependent external force $g\in W^{1,2}_{\text{loc}}(\mathbb R; L^2(\mathbb R))$ satisfies
\begin{equation}
\|g(t)\|_{L^2(\mathbb R)}^2 + \|g'(t)\|_{L^2(\mathbb R)}^2 \leq
Me^{\gamma |t|}, \quad t\in\mathbb R,
\label{e2}
\end{equation}
with $\gamma <\lambda$ and $M$ is independent of $t$. Note that, when
$\mu  = 0$, equation \eqref{e1} becomes the well-known weakly damped forced KdV
equation.

\medskip
The KdV-Burgers equation has been derived as an asymptotic model for
the propagation of weakly nonlinear dispersive long waves in certain
physical situations when dissipative and damping effects occur, such as waves in
an elastic tube, ion-sound waves damped by ion-neutron collisions,
waves in a tube filled with viscous fluid, flow of a fluid containing
gas bubbles, turbulence, or crystal lattice theory (see e.g.
\cite{Joh70, Joh72, OttSudan, SuGadner}) and as a model for
propagation of tsunamis (see \cite{Dutykh}). The well-posedness of homogeneous KdV-Burger equations, i.e. with $g \equiv 0$, has been thoroughly studied (see e.g. \cite{molinet2011sharp,molinet2013sharp,cavalcanti2014global} and references therein). The large time dynamics of KdV-Burgers in the phase spaces $H^1(\mathbb R)$ or $H^2(\mathbb R)$ have also been extensively investigated, see e.g. \cite{Dlotko,DlotkoSun} for the autonomous case and \cite{AB12} for the non-autonomous case. Concerning the phase space $L^2(\mathbb R)$ and autonomous external force $g \equiv g(x)$, the existence, regularity and finite fractal dimension of a global attractor are obtained in \cite{boling1998global}, \cite{guo2017regular} and \cite{guo2020dimension}, respectively. It is also worth mentioning the related works with low regular initial data such as \cite{goyal2022global,mcconnell2023global}.  In practice, the external force is likely to be time dependent, i.e. $g \equiv g(t,x)$, leading to the non-autonomous equation \eqref{e1}. The long-time behavior of solutions this equation in phase space $L^2(\mathbb{R})$ has not been studied, and this serves as a motivation for the current paper.

\medskip
The first objective of this paper is to study the existence and
regularity of a pullback attractor in 
 %$L^2(\mathbb{R})$
 for the non-autonomous
process associated to problem \eqref{e1}.
% in the case that the external
%force $g$ depends on the time $t$. 
More precisely, we show that the pullback attractor
of the associated process exists in $L^2(\mathbb{R})$ and is in fact
compact in $H^3(\mathbb{R})$; thus, obtaining an asymptotic smoothing
effect of solutions which generally belong only to
$L^2(\mathbb{R})$. Note that the result obtained are of optimal order
and holds for both the KdV-Burgers equation (i.e., when $\mu >0$) and
the weakly damped KdV equation (i.e., when $\mu =0$); moreover, when
$\mu =0$ and $g\equiv g(x) \in L^2(\mathbb{R})$, we recover the results of
\cite{GoR} for the autonomous weakly damped forced KdV equation. 
Since the initial data $u_\tau$ is in
$L^2(\mathbb R)$, studying well-posedness and estimates of the solutions requires Bourgain-space
techniques. Furthermore, the problem under consideration is posed in the real line $\mathbb R$, which makes Sobolev embeddings $H^1(\mathbb R)\hookrightarrow L^2(\mathbb R)$ no longer
compact. These lead to certain difficulties for proving the pullback
asymptotic compactness of solutions and therefore in studying the
existence and regularity of pullback attractors for the KdV-Burgers
equation. In this paper, we deal with these issues by first using the
energy equation method introduced by Ball in \cite{Ball} to get the
existence of a family of bounded pullback absorbing sets in
$L^2(\mathbb{R})$. To obtain the higher regularity of the pullback
attractor, following the general lines of the
approach in \cite{Goubet,GoR}, we split the solutions into two parts possessing different frequencies,
in which the first one decays exponentially in $L^2(\mathbb R)$ while the other one is
bounded in $H^3(\mathbb R)$. Once the boundedness in
$H^3(\mathbb{R})$ of the regular part of the solution is obtained, the compactness
of the pullback attractor in $H^3(\mathbb{R})$ the follows by applying again
the energy equation method, now for the equation 
$u_t$ along trajectories lying in the pullback attractor. Delicate issues
arise due to the lack of compact Sobolev embeddings which are
compensated by exploiting the dispersive regularization property of
the equation under study in Bourgain spaces.

\medskip
The second objective of this paper is to study the upper
semi-continuity of the pullback attractors  $\{\mathcal
A_{\mu }\}_{\mu \in[0,1]}$ at $\mu  = 0$ and their finite fractal dimensions in $H^1(\mathbb R)$.  
  In other words, we show that when $\mu  \rightarrow 0^+$, the pullback attractor $\mathcal
A_\mu $ of the KdV-Burgers equation converges to the pullback
attractor~$\mathcal A_0$ of the weakly damped KdV equation in the
sense of the Hausdorff semi-distance in $L^2(\mathbb R)$. This result
is new even in the autonomous case. This is done by showing the compactness of the pullback attractors \textit{uniformly in $\mu\ge 0$} using some tail estimates technique. Finally, we show that the pullback attractor has a finite fractal dimension in $H^1(\mathbb R)$ for any $\mu > 0$. In order to do this, we prove a result for finite fractal dimensions for non-autonomous problem, which is an extension of the autonomous case in \cite[Theorem 4.15]{CL08}. Since this extension is well applicable for \eqref{KdV} on the whole line $\R$ and avoids using finite dimensional projection, we expect that it can be used in other non-autonomous problems in unbounded domains. Furthermore, by exploiting the dissipation by the term $-\mu u_{xx}$, we can get our result without having to pursue higher regularity of the pullback attractors, as in \cite{WH20}, which requires to work with associated Bourgain spaces and likely imposes higher regularity of the external force. 

\medskip
As a final remark, we believe that the results of this paper can be extended to show that if $g(t)\in
H^k(\mathbb{R})$, for some $k\in \mathbb{N}$, then equation \eqref{e1}
generates a process in the phase space $H^k(\mathbb{R})$ which has a
pullback attractor in $H^k(\mathbb{R})$, and the pullback attractor is
compact in $H^{k+3}(\mathbb{R})$. Therefore, this paper, in a sense, completes the study of
existence and regularity of pullback attractors for equation
\eqref{e1} in the scale of Sobolev spaces $H^k(\mathbb{R})$ for
integers $k\geq 0$. The details are left for the interested reader.

\medskip
\textbf{The paper is organized as follows}:  In section \ref{sec:2} show the global existence of \eqref{e1} in the phase space $L^2(\mathbb R)$. Section \ref{sec:KdV-Burgers} is devoted to the the existence and regularity of pullback attractor $\mathcal A_{\mu}$ for the for the process generated by \eqref{e1}. In the last two sections, we prove the upper semicontinuity of $\mathcal A_{\mu}$ at $\mu = 0$ and their finite fractal dimensions, respectively.

%\vskip 2cm
\medskip
Throughout this paper, we use the following notation:
\begin{itemize}
	\item The usual norm in $L^2(\mathbb R)$ is written as $\|\cdot\|$ or convenience.
	\item $\mathcal S(\mathbb R)$ denotes the Schwartz space on $\mathbb R$ and $\mathcal S'(\mathbb R)$ is the corresponding space of distributions. For $f\in \mathcal S'$,
	$\widehat{f}$ or $\mathcal F(f)$ is its Fourier transform, i.e.,
	\begin{equation*}
		\widehat{f}(\xi) = \int_{\mathbb R^n}e^{-i z\cdot \xi }f(z)\,dz.
		\label{e3}
	\end{equation*}
	The usual Sobolev space is denoted by $H^{s}(\mathbb R)$ for
	$s\in\mathbb R$ with norm
	\begin{equation*}
		\|u\|_{H^{s}} = \|\langle\xi\rangle^{s}\widehat{f}(\xi)\|_{L^2(\mathbb R)},
		\label{e4}
	\end{equation*}
	where $\langle\xi\rangle = \sqrt{1+|\xi|^2}$. 
	
	\item The Bourgain space-time
	function space $X^{s,b}$ for $s,b\in \mathbb{R}$ is defined as the
	completion of the Schwartz space $\mathcal{S}(\mathbb{R}^2)$ with
	respect to the norm
	\begin{equation}\label{BourgainSpace}
	\begin{aligned}
		\|u\|_{X^{s,b}} &= \left\|\langle i(\tau-\xi^3) + \xi^2\rangle^{b}\langle\xi\rangle^{s}
		\widehat{u}(\tau,\xi)\right\|_{L^{2}(\mathbb R^2)}\\
		&=\left\|\langle |\tau-\xi^3| + |\xi^2|\rangle^{b}\langle\xi\rangle^{s}
		\widehat{u}(\tau,\xi)\right\|_{L^{2}(\mathbb R^2)}.
	\end{aligned}
	\end{equation}

	We have the continuous embedding
	\begin{equation*} \label{e5a}
		X^{s,b} \subset C_b(\mathbb{R}; H^s(\mathbb{R}_x))
	\end{equation*}
	for any $s\in \mathbb{R}$ and $b>{1}/{2}$. For an
	interval $I\subseteq \mathbb R$, we define the space $X_{I}^{s,b}$
	with norm
	\begin{equation*}
		\|u\|_{X_{I}^{s,b}}=\inf_{w\in X^{s,b}}\{\|w\|_{X^{s,b}}\colon 
		w \equiv u \text{ on } I\}.
		\label{e7}
	\end{equation*}

	\item For a Banach space $(X, \|\cdot\|)$, $\mathcal{B}(X)$ denotes the set of
	bounded subsets of $X$, and for $A,\,B\subseteq X$, the Hausdorff
	semi-distance between $A$ and $B$ is defined by
	\[
	\operatorname{dist}_{X}(A,B)=\sup_{x\in A}\inf_{y\in B}\|x-y\|.
	\] 
\end{itemize}

\section{Global existence of solutions}\label{sec:2}

We start with some results on the local well-posedness of the KdV-Burgers
equation in spaces of low regularity. The results presented here follow from \cite{MoRi02}.
Denote by $W$ the group associated with the free evolution
of \eqref{e1}, i.e. $u_t + u_{xxx} - \mu  u_{xx}= 0$, 
\begin{equation*}
	\mathcal F(W(t)\varphi)(\xi) = \exp(-\mu  |\xi|^2t +
	i\xi^3t)\widehat{\varphi}(\xi), \;\;\varphi \in \mathcal S'(\mathbb R), \;\;\forall
	t\geq 0,
	\label{e8}
\end{equation*}
and we extend $W$ to the whole real axis by
\begin{equation*}
	\mathcal F(W(t)\varphi)(\xi) = \exp(-\mu  |\xi|^2|t| +
	i\xi^3t)\widehat{\varphi}(\xi), \;\;\varphi \in \mathcal S'(\mathbb R), \;\;\forall
	t\in\mathbb R.
	\label{e9}
\end{equation*}
Fix $0 \leq s \leq 23/48$ and $1/2 < b < 25/48$. Assume that $\lambda >0$ and $g = g(x,t) \in X^{-s,b-1}_{[\tau - T_g, \tau + T_g]}$ for some $T_g>0$. We look for a
fixed-point of the mapping 
\begin{equation}
	\begin{aligned}
		\mathcal{M}(u)(t) = \psi_{T_g}(t)\biggl(W(t)u_{\tau} - \frac{\chi_{\mathbb
				R_{\tau}}(t)}{2}&\int_{\tau}^{t}W(t-s)\biggl[(\psi_{T_g}(s)u(s))^2_{x}\\
		& \quad + 2\lambda \psi_{T_g}(s)u(s) - 2\psi_{T_{g}}(s)g \biggr]ds\biggr)
		%    
		%    \left(
		%W(t-s)\partial_{x}(\psi_{T}^2(s)u^2(s))\right. \\ & \left. +\,\mu 
		%\psi_{T}^2(s)u_{xx}(s) - \lambda \psi_{T}^2(s)u(s) -
		%\psi_{T}^2(s)g(s)\right)ds\biggr),
	\end{aligned}
	\label{k1}
\end{equation} 
where $\mathbb R_{\tau} = [\tau,+\infty)$ and $\chi_{\mathbb
	R_{\tau}}$ is its characteristic function, $\psi$ is a time
cut-off function satisfying
\begin{equation*}
	\psi\in C_{0}^{\infty}(\mathbb R),\; \operatorname{supp}\psi \subseteq [\tau -1,\tau +1],\;
	\psi \equiv 1 \text{ on } \left[\tau - \frac{1}{2},\tau + \frac{1}{2}\right],
\end{equation*}
and $\psi_{T}(\cdot) = \psi(\cdot/T)$. Indeed, if $u$ solves
\eqref{k1} on $[\tau-T,\tau+T]$, then $u$ is a solution of \eqref{e1}
on $[\tau,\tau+T/2]$. %The following result was proved in 

\begin{theorem}[\cite{GoR}]\label{solution}
	Let $0\leq s\leq 23/48$ and $1/2< b<25/48$. For any $u_\tau \in
	L^2(\mathbb R)$, there exists $T>0$ such that the problem \eqref{e1}
	is uniquely solvable in $X_{[\tau,\tau+T]}^{-s,b}$. In particulars, the solution $u\in C([\tau,T];L^2(\mathbb R))$. Moreover, this solution can be extended globally on $\mathbb R_\tau$.
	%Moreover, the map $u_\tau \mapsto u$ is analytic from $L^2(\mathbb R)$ to  $X_{[\tau,\tau+T]}^{-s,b}$.
\end{theorem}
\begin{remark}
	It is in fact possible to show sharper local well-posedness results of \eqref{e1}, e.g. \eqref{e1} is well-posed in $H^{s}(\mathbb R)$ for any $s>-1$. The interested reader is referred to \cite{MoRi02} for more details.
	%In the present paper, however, we use the Bourgain space $X^{s,b}$ as
	%defined in \eqref{e5}, as we want to apply the bilinear estimate
	%arguments of \cite{GoR} to prove the decay of the part $w$ of the
	%solution $u$ (for details see Section~3.1.1).
\end{remark}
\begin{proof}
	The proof of local well-posedness follows from \cite{MoRi02} with slight modifications, so we omit it here. To show the global existence, we take the inner product of \eqref{e1} with $2u(t)$ in $L^2(\mathbb R)$ to get
	\begin{equation}
		\frac{d}{dt}\|u(t)\|^2 + 2\lambda\|u(t)\|^2+2\mu \|u_x(t)\|^2 = 2(g(t), u(t)).
		\label{f1}
	\end{equation}
	Applying the Cauchy and Young's inequalities, we see that
	\begin{equation*}
		\frac{d}{dt}\|u(t)\|^2 + \frac{3\lambda}{2}\|u(t)\|^2 + 2\mu \|u_x(t)\|^2 \leq \frac{1}{2\lambda}\|g(t)\|^2.
		%\label{f2}
	\end{equation*}
	Therefore,
	\begin{equation}
		\begin{aligned}
			\|u(t)\|^2 &+ \int_{\tau}^{t}e^{-\lambda (t-s)}\left(\frac{\lambda}{2}\|u(s)\|^2 + 2\mu \|u_x(s)\|^2\right)ds\\
			&\leq e^{-\lambda(t-\tau)}\|u_{\tau}\|^2 + \frac{e^{-\lambda t}}{2\lambda}\int_{\tau}^{t}e^{\lambda s}\|g(s)\|^2ds.
		\end{aligned}
		\label{f3}
	\end{equation}
	This implies that the local solution can be extended globally on $\mathbb R_{\tau}$.
\end{proof}

\section{Existence and regularity of pullback attractors}\label{sec:KdV-Burgers}
This section is devoted to the existence and regularity of pullback attractors for KdV-Burgers equation, that is \eqref{e1} with $\mu >0$. We first recall an useful bilinear estimate which is crucial in proving the regularity of the attractors. Next, the existence of the pullback attractor in $L^2(\mathbb R)$ is shown in subsection \ref{L2attractor} by an energy equation estimate. The boundedness and compactness of the attractor in $H^3(\mathbb R)$ are proved in subsections \ref{boundednessH3} and \ref{regularity} respectively.

%\begin{remark}\rm
%From Theorem~\ref{solution} and \eqref{e5a}, we infer that problem
%\eqref{e1} is locally well-posed in $H^{s}(\mathbb R)$ for $-23/48\leq
%s\leq 0$. In particular, this problem is locally well-posed in
%$L^2(\mathbb R)$.  This is the space in which we shall prove the
%global well-posedness and also the existence of a pullback
%attractor. The well-posedness of \eqref{e1} in $X^{s,b}$ for
%$-23/48\leq s\leq 0$ is then used in the proof of
%Theorem~\ref{compactness}, where we establish the compactness of the
%pullback attractor in $H^3(\mathbb R)$.
%\end{remark}

\medskip
We first briefly recall the basic concepts concerning pullback attractors fo non-autonomous dynamical systems. For more details, the interested reader is referred to the book \cite{CLR13}.

Let $\{U(t,\tau)\colon t\geq \tau,\, \tau\in \mathbb{R}\}$ be a
process in $X$, i.e., $U(t,\tau)\colon X\rightarrow X$ satisfy
$U(\tau,\tau)=\operatorname{Id}_X$ and $U(t,s)U(s,\tau)=U(t,\tau)$ for
all $t\geq s\geq \tau$, $\tau \in \mathbb{R}$. This process is called \textit{continuous} (\textit{weakly continuous}) if for each $t>\tau$ the mapping $x \mapsto U(t,\tau)x$ is continuous (weakly continuous). 
\begin{definition}\label{defpullbackattractor}
	A family $\mathcal A=\{A(t)\colon t\in \mathbb{R}\}\subset
	\mathcal{B}(X)$ is said to be a \textit{pullback attractor} 
	for the process $\{U(t,\tau)\}$ if
	\begin{itemize}
		\item[(1)] $A(t)$ is compact for all $t\in \mathbb{R}$, 
		\item[(2)] $\mathcal A$ is invariant, i.e., 
		\[
		U(t,\tau)A(\tau)=A(t) \text{ for all } t\geq \tau.
		\]
		\item[(3)] $\mathcal A$ is {\it pullback attracting}, i.e., 
		\[ 
		\lim\limits_{\tau\rightarrow
			-\infty}\operatorname{dist}_{X}(U(t,\tau)D,A(t))=0 \text{ for all }
		D\in \mathcal{B}(X) \text{ and all } t\in \mathbb{R}.
		\]
		\item[(4)] If $\{C(t)\colon t\in \mathbb{R}\}$ is another family
		of closed pullback attracting sets, then $A(t)\subset C(t)$
		for all $t\in \mathbb{R}$.
	\end{itemize}
\end{definition}
\begin{theorem}[\cite{LiZhong, CLR13}]\label{main}
	Let $\{U(t,\tau)\}$ be a process satisfying the following
	conditions\/\textup{:}
	\begin{itemize}
		\item[(i)] $\{U(t,\tau)\}$ is continuous or weakly continuous;
		\item[(ii)] there exists a family of pullback absorbing sets
		$\mathcal{B} = \{B(t)\colon t\in \mathbb{R}\}$,
		i.e., for any $t\in \mathbb{R}$ and any $D\in \mathcal{B}(X)$, there
		is a $\tau_0 = \tau_0(D,t)\leq t$ such that
		\[
		\bigcup\limits_{\tau \leq \tau_0}U(t,\tau)D \subseteq B(t);
		\]
		\item[(iii)] $\{U(t,\tau)\}$ is pullback asymptotically compact, i.e. for any $t\in \mathbb{R}$, any
		$D\in \mathcal{B}(X)$, any sequence $\tau_n\rightarrow -\infty$, and
		any sequence $\{x_n\}\subset D$, the sequence $\{U(t,\tau_n)x_n\}$ is
		relatively compact.
	\end{itemize}
	Then $\{U(t,\tau)\}$ has a unique pullback attractor $\mathcal A =
	\{A(t)\colon t\in \mathbb{R}\}$, and
	\[
	A(t) = \bigcap_{s\leq t}\overline{\bigcup_{\tau \leq s}U(t,\tau)B(\tau)}^{\,X},
	\]
	where $\bar{A}^{\,X}$ denotes the closure of $A$ with respect to
	the norm topology of $X$.
\end{theorem}

\subsection{Pullback attractors in $L^2(\mathbb R)$}\label{L2attractor}

Thanks to Theorem \ref{solution}, for each $\mu \ge 0$, we can define a process $\{U_\mu(t,\tau)\}$ in $L^2(\mathbb R)$ as follows: for any $u_\tau \in L^2(\mathbb R)$, $U(t,\tau)u_\tau$ is the solution to \eqref{e1} established in Theorem \ref{solution} at time $t>\tau$ with initial data $u_\tau$.

\begin{proposition}\label{absorbingset}
There exists a family of pullback absorbing sets $\mathcal B = \{B(t): t\in \mathbb R\}$ in $L^2(\mathbb R)$ for $U_{\mu }(t, \tau)$.
\label{absorbing}
\end{proposition}
\begin{proof}

Denote by $B(t)$ the ball in $L^2(\mathbb R)$ centered at origin and has the radius
\begin{equation}\label{r}
	r(t) = \left(\frac{e^{-\lambda
	    t}}{\lambda}\int_{-\infty}^{t}e^{\lambda
	  s}\|g(s)\|^2ds\right)^{1/2}.
\end{equation}
Then the estimate \eqref{f3} implies that 
for any bounded set $B\subset L^2(\mathbb
R)$ there exists $\tau_0\in \mathbb R$ such that
\begin{equation*}
U_{\mu }(t,\tau)u_{\tau}\in B(t)  \qquad \forall \tau\leq\tau_0,\quad \forall u_\tau
  \in B.
%\label{f4}
\end{equation*}
\end{proof}

We prove in the next lemma the weak continuity of the process
$\{U_{\mu }(t,\tau)\}$, which is not only needed for the abstract Theorem \ref{main} but also is very useful in showing the asymptotic compactness via an energy equation technique.
\begin{lemma}\label{weakcontinuous}
The process $\{U_{\mu }(t,\tau)\}$ is weakly continuous on
$L^2(\mathbb{R})$ in the sense that, if $u^{\epsilon}_{\tau}
\rightharpoonup u_{\tau}$ in $L^2(\mathbb R)$, then
$U_{\mu }(t,\tau)u^{\epsilon}_{\tau} \rightharpoonup
U_{\mu }(t,\tau)u_{\tau}$ in $L^2(\mathbb R)$ for all $t\geq \tau$.
\end{lemma} 
\begin{proof}
The proof goes similarly to \cite[Section 6]{Goubet}. Without loss of generality, we consider the case $\lambda = 0$, $g =
0$, and $\mu  = 1$. We fix $N$ and define the projection $P_N$ by
\begin{equation}\label{PN}
\mathcal F(P_Nu)(\xi) = \widehat{u}(\xi)\chi\left(\frac{\xi}{N}\right),
\end{equation}
where $\chi$ is the characteristic function of $[-1, 1]$. Let $v^{\epsilon}_{N}(t)$ be the solution of 
\begin{equation}
\partial_tv^{\epsilon}_{N} + \partial_{xxx}v^{\epsilon}_{N} + v^{\epsilon}_{N}\partial_xv^{\epsilon}_{N} -\partial_{xx}v^{\epsilon}_{N} = 0,
\label{f5} 
\end{equation}
with initial data $v^{\epsilon}_{N}(\tau) =
P_Nu^{\epsilon}_{\tau}$. Since $P_Nu^{\epsilon}_{\tau}$ is smooth and
\begin{equation*} 
\|P_Nu^{\epsilon}_{\tau}\|_{H^{2}(\mathbb R)} \leq KN^2,
%\label{f6}
\end{equation*}
where $K$ is independent of $\epsilon$, we conclude that $v^{\epsilon}_N$ is bounded in 
\begin{equation*}
C\left([\tau-T,\tau+T];H^2(\mathbb R)\right)\cap C^{1}\left([\tau-T,\tau+T];H^{-1}(\mathbb R)\right),
%\label{f7}
\end{equation*}
as a result of \cite{temam1969probleme}. Thus, 
\begin{equation*}
v^{\epsilon}_{N}\rightarrow v_N \text{ strongly in } C\left([\tau-T,\tau+T];L^2_{\text{loc}}(\mathbb R)\right),
%\label{f8}
\end{equation*}
because $H^2(\mathbb R) \hookrightarrow L^2_{\text{loc}}(\mathbb R)$
with a compact embedding. Hence, in equation \eqref{f5} we can pass to
the limit to find that $v_N$ is a solution of \eqref{f5} with initial
data $v_N(\tau) = P_Nu_{\tau}$. On the other hand,
\begin{equation*}
v_N\in L^{\infty}([\tau-T,\tau+T];H^2(\mathbb R)) \text{ and }
\partial_t v_N \in L^{\infty}([\tau-T,\tau+T];H^{-1}(\mathbb R)).
%\label{f9}
\end{equation*}
We conclude that $v_{N}\in C([\tau-T,\tau+T];H^{s}(\mathbb
R))$ for all $s<2$.

Now, let $w^{\epsilon}_{N}(t) = u^{\epsilon}(t) - v^{\epsilon}_N(t)$,
where $u^{\epsilon}(t) = U_{\mu }(t,\tau)u^{\epsilon}_\tau$. We get
\begin{equation}
\partial_tw^{\epsilon}_{N} + \partial_{xxx}w^{\epsilon}_{N} +
w^{\epsilon}_{N}\partial_{x}w^{\epsilon}_{N} +
(v^{\epsilon}_{N}w^{\epsilon}_{N})_{x} - \partial_{xx}w^{\epsilon}_{N}
= 0,
\label{f10}
\end{equation}
with initial data $w^{\epsilon}_{N}(\tau) = (\operatorname{Id} -
P_N)u^{\epsilon}_{\tau}$ (including the case $\epsilon = 0$). From the
well-posedness of \eqref{f10} in $X^{-s,0}$, $s\in [0,3/4)$, we obtain
that on bounded time intervals,
\begin{equation}
\|w^{\epsilon}_{N}(t)\|_{H^{-s}(\mathbb R)} \leq K\|w^{\epsilon}_{N}(\tau)\|_{H^{-s}(\mathbb R)},
\label{f11}
\end{equation}
where $K$ depends on $T$, but is independent of $N$ and $\epsilon$. On
the other hand, since $\{u^{\epsilon}_{\tau}\}$ is bounded in
$L^2(\mathbb R)$,
\begin{equation*}
\|w^{\epsilon}_{N}(\tau)\|_{H^{-s}(\mathbb R)} = \|(\operatorname{Id}
- P_N)u^{\epsilon}_{\tau}\|_{H^{-s}(\mathbb R)} \leq
\frac{\|u^{\epsilon}_{\tau}\|_{L^{2}(\mathbb R)}}{N^{s}} \leq KN^{-s}.
%\label{f12}
\end{equation*}
Now, for a test function $\psi$ that is smooth, compactly supported,
and satisfies $\|\psi\|_{L^2(\mathbb R)} = 1$,
\begin{equation*}
\begin{aligned}
\left|(u(t)-u^{\epsilon}(t),\psi)\right| &= \left|\left(w^{0}_{N}(t) + v_N(t)-v^{\epsilon}_N(t) + w^{\epsilon}_N(t), \psi\right)\right|\\
&\leq 2KN^{-s} + \left|(v_N(t) - v^{\epsilon}_N(t), \psi)\right|
\end{aligned}
%\label{f13}
\end{equation*}
because of \eqref{f11}. From the estimates on $\{v^{\epsilon}_{N}\}$,
letting $\epsilon\rightarrow 0$ with $N$ fixed we obtain
\begin{equation*}
\limsup_{\epsilon\rightarrow 0}\left|(v_N(t)-v^{\epsilon}_N(t),\psi)\right| \leq 2KN^{-s}.
%\label{f14}
\end{equation*}
Let $N\rightarrow +\infty$ to conclude that $u^{\epsilon}(t)$
converges to $u(t)$ in the distributional sense.  Since
$w^{\epsilon}_N$ is a bounded sequence in $L^{\infty}([\tau-T,\tau+T];
L^2(\mathbb R))$, a density argument completes the proof.
\end{proof}
\begin{theorem}\label{L2}
The process $\{U_{\mu }(t,\tau)\}$ possesses a pullback attractor
$\mathcal A_{\mu } = \{A_{\mu }(t)\colon$ \linebreak $t\in \mathbb R\}$ in
$L^2(\mathbb R)$.
\end{theorem}
\begin{proof}
By Lemma \ref{weakcontinuous} and Proposition \ref{absorbingset}, $\{U_{\mu }(t,\tau)\}$ is weakly continuous and possesses a
family of bounded pullback absorbing set $\mathcal B = \{B(t)\colon
t\in \mathbb R\}$. Thanks to Theorem \ref{main}, it remains to show the pullback asymptotic
compactness of $\{U_{\mu }(t,\tau)\}$. Let $\{x_n\}$ be a bounded
sequence in $L^2(\mathbb{R})$ and $\{\tau_n\} \subset \mathbb R$ be a sequence such
that $\tau_n \rightarrow -\infty$. We will prove that
$\{U_{\mu }(t,\tau_n)x_n\}$ is precompact in $L^2(\mathbb R)$.

From Proposition \ref{absorbingset}, we have that
$\{U_{\mu }(t,\tau_n)x_n\}$ is bounded in $L^2(\mathbb R)$. Thus there
exists a $w_0$ belonging to $L^2(\mathbb R)$ such that
\begin{equation*}
U_{\mu }(t,\tau_n)x_n \rightharpoonup w_0 \text{ as } n \rightarrow +\infty.
%\label{f15}
\end{equation*}
For fixed $T>0$, by a similar argument, we get a $w_T\in L^2(\mathbb
R)$ such that
\begin{equation*}
U_{\mu }(t-T,\tau_n)x_n \rightharpoonup w_T \text{ as } n \rightarrow +\infty.
%\label{f16}
\end{equation*}
Then, because of the weak continuity of $\{U_{\mu }(t,\tau)\}$, for any $T>0$,
\begin{equation*}
\begin{aligned}
w_0&= \wlim U_{\mu }(t,\tau_n)x_n= \wlim U_{\mu }(t,t-T)U_{\mu }(t-T,\tau_n)x_n\\
&= U_{\mu }(t,t-T)\wlim U_{\mu }(t-T,\tau_n)x_n= U_{\mu }(t,t-T)w_T,
\end{aligned}
%\label{f17}
\end{equation*}
where $\wlim$ denotes the weak limit in $L^2(\mathbb R)$ as $n\rightarrow +\infty$. From
\eqref{f1}, denoting $E(u) = \lambda \|u\|^2 + 2\mu \|u_x\|^2$, we get
\begin{multline}
\|U_{\mu }(t,\tau)u_{\tau}\|^2 + \int_{\tau}^te^{-\lambda (t-s)}E(u)ds 
= e^{-\lambda(t-\tau)}\|u_{\tau}\|^2 + 2\int_{\tau}^{t}e^{-\lambda(t-s)}(g(s), u(s))\, ds.
\label{f18}
\end{multline}
Applying \eqref{f18} to $w_0 = U_{\mu }(t,t-T)w_T$ and
$U_{\mu }(t,\tau_n)x_n = U_{\mu }(t,t-T)U_{\mu }(t-T,\tau_n)x_n$
respectively, we get
\begin{equation}
\begin{aligned}
\|w_0\|^2 &+ \int_{t-T}^te^{-\lambda(t-s)}E(U_{\mu }(s,t-T)w_T)\,ds\\
&= e^{-\lambda T}\|w_T\|^2 + 2\int_{t-T}^te^{-\lambda(t-s)}
\int_{\mathbb R}(g(s),U_{\mu }(s,t-T)w_T)\,ds
\end{aligned}
\label{f19}
\end{equation}
and
\begin{equation*}
\begin{aligned}
\|U_{\mu }(t,\tau_n)x_n\|^2 +
&\int_{t-T}^{t}e^{-\lambda(t-s)}E(U_{\mu }(s,\tau_n)x_n)\,ds\\ &=
e^{-\lambda T}\|U_{\mu }(t-T,\tau_n)x_n\|^2 +
2\int_{t-T}^{t}e^{-\lambda(t-s)}(g(s), U_{\mu }(s,\tau_n)x_n)\,ds.
%\label{f20}
\end{aligned}
\end{equation*}
It is easy to see that 
\begin{equation*}
u \mapsto \left(\int_{t-T}^{t}e^{-\lambda(t-s)}E(u)ds\right)^{1/2}
%\label{f21}
\end{equation*}
is a norm in $L^{2}(t-T,t;H^1(\mathbb R))$ that is equivalent to the
usual one. Then, from \eqref{f3}, we get as $n\rightarrow +\infty$
\begin{equation*}
	U_{\mu }(\cdot, \tau_n)x_n \rightharpoonup U_{\mu }(\cdot, t - T)w_T \quad \text{ weakly in } L^2(t-T, t;H^1(\mathbb R))
\end{equation*}
and consequently
\begin{equation}
\int_{t-T}^{t}e^{-\lambda(t-s)}E(U_{\mu }(\cdot, t-T)w_T)\,ds \leq 
\limsup_{n\rightarrow \infty}\int_{t-T}^{t}e^{-\lambda(t-s)}E(U_{\mu }(\cdot, \tau_n)x_n)\,ds.
\label{f22}
\end{equation}
From \eqref{f19} and \eqref{f22}, and noting that
$U_{\mu }(\cdot,\tau_n)x_n \rightharpoonup U_{\mu }(\cdot, t-T)w_T$ in
$L^{2}(t-T,t;L^2(\mathbb R))$, we get that
\begin{equation*}
\begin{aligned}
\|w_0\|^2 &\geq \limsup_{n\rightarrow\infty}\|U_{\mu }(t,\tau_n)x_n\|^2
- \limsup_{n\rightarrow\infty}e^{-\lambda
  T}\|U_{\mu }(t-T,\tau_n)x_n\|^2\\ &\geq \limsup_{n\rightarrow
  \infty}\|U_{\mu }(t,\tau_n)x_n\|^2 + e^{-\lambda T}r(t-T)^2
\end{aligned}
%\label{f23}
\end{equation*}
where $r(t)$ is defined as in \eqref{r}. Hence, thanks to $\lim_{T\rightarrow+\infty}e^{-\lambda T}r(t-T)^2 = 0$, we obtain
\begin{equation*}
\|w_0\|^2 \geq \limsup_{n\rightarrow \infty}\|U_{\mu }(t,\tau_n)x_n\|^2.
%\label{f24}
\end{equation*}
Combined with the fact that $U_{\mu }(t,\tau_n)x_n \rightharpoonup w_0$
in $L^2(\mathbb R)$, we have
\begin{equation*}
U_{\mu }(t,\tau_n)x_n \rightarrow w_0 \text{ strongly in } L^2(\mathbb
R),
%\label{f25}
\end{equation*}
when passing to a subsequence if necessary. This completes the proof.
\end{proof}

\subsection{Boundedness in $H^{3}(\mathbb R)$}\label{boundednessH3}
In this subsection, we will prove that the $L^2(\mathbb R)$ pullback attractor $\mathcal A_{\mu }$ of \eqref{e1} obtained in subsection \ref{L2attractor} is in fact bounded in $H^3(\mathbb R)$. In order to do that, we will employ the method of splitting solution. Roughly speaking, we will divide the solution into two parts in which one part decays exponentially in $L^2(\mathbb R)$ while the other part is bounded in $H^3(\mathbb R)$.

\subsubsection{Splitting solutions}
We write $u = v+w$, where $v, w$ are solutions of the equations
\begin{equation}
	v_t  + v_{xxx} - \mu  v_{xx}+ \lambda v + vv_x = g(t) - P_N((vw)_{x} + ww_{x}),\quad x\in\mathbb R,\; t>\tau
	\label{g1}
\end{equation}
\begin{equation}
	w_t + w_{xxx} - \mu  w_{xx} + \lambda w  + Q_N(ww_x) = -Q_N(vw)_{x}, \quad x\in\mathbb R,\; t>\tau
	\label{g2}
\end{equation}
with initial data $v(\tau) = P_Nu_{\tau}$ and $w(\tau) =
Q_{N}u_{\tau}$, where we recall that $P_N$ is the projection defined in \eqref{PN} and $Q_N = Id - P_N$.

Since \eqref{g2} is a dissipative equation without forces, its solution $w$ is expected to be exponentially decay. While the solution to \eqref{g1} is expected to be smooth, i.e. belong to $H^3(\mathbb R)$, thanks to the regularity of initial data $v(\tau) = P_Nu_\tau$.
%First, we prove the global existence of the part $w$ and the decay in
%time of $w(t)$ in $L^2(\mathbb R)$. That implies the global existence
%of $v$. Then we prove boundedness of $v$ in $H^3(\mathbb R)$.

The analysis of the rest of this section frequently uses the following bilinear estimate.
\begin{lemma}\label{bilinear2}
	Let $s\in \left(-1/2, 0\right]$ be given. Then there exists a
	constant $c_{s}$ such that, for $b' \in (1/2, 9/16]$, $b\in [1/2,
	\min\{b', 1+s\}]$, and any function $u$ in $X^{s,b}$,
	\begin{equation}
		\|D(u^2)\|_{X^{s,b'-1}} \leq c_{s}\|u\|_{X^{s,b}}\|u\|_{X^{-1/2,b}}.
		\label{est1}
	\end{equation}
	
	As a corollary, if $u \in X^{s,b}$ is a function satisfying that the Fourier transform $\widehat{u}(\xi, \tau)$
	of $u$ is supported in $\{(\xi, \tau)\colon |\xi|\geq N\}$, where
	$N>0$. Then, for $1/2<b\leq 9/16$,
	\begin{equation}
		\|D(u^2)\|_{X^{0,b-1}}\leq \frac{C}{\sqrt{N}}\,\|u\|_{X^{0,b}}^2.
		\label{est2}
	\end{equation}
\end{lemma}
\begin{proof}
	The estimate \eqref{est1} is similar to \cite[Proposition 2.2]{Goubet} for the KdV equation. However, we emphasize that we currently use different Bourgain spaces, which are based on the symbol of the free KdV-Burgers equation instead of the free KdV equation. The proof of \eqref{est1} then follows instead similarly from \cite[Proposition 3.2]{MoRi02} with slight modifications.
	%
	%Lemma \ref{bilinear2} needs a new proof. We postpone a proof of Lemma \ref{bilinear2} to the Appendix \ref{bilinearproof} to continue with the regularity of the pullback attractor.
\end{proof}

\subsubsection{Well-posedness and decay of the part $w$}
Since $v = u - w$, we can rewrite \eqref{g2} as
\begin{equation}
	w_t + w_{xxx} + \lambda w - \mu  w_{xx} - Q_N(ww_x) = -Q_N(uw)_{x}
	\label{g3}
\end{equation}
with $w(\tau) = Q_Nu_{\tau}$. The local well-posedness of \eqref{g3}
follows as for \eqref{e1}. Applying arguments of Section~2 to the
equation
\begin{equation*}
	\begin{aligned}
		w(t) &= \psi_1(t)W(t)w(\tau) \\
		&\quad + \frac{1}{2}\,\psi_{T}(t)
		\int_{\tau}^{t}W(t-s)\biggl[-2\lambda \psi_1(s)w(s)- \left(Q_N(\psi_{T}(s)u(s)\psi_{T}(s)w(s))\right)_{x} \\
		&\quad \hspace{4cm}
		+ \left(Q_N(\psi_T(s)w(s))^2\right)_{x}\biggr]ds.
	\end{aligned} 
	%\label{g4}
\end{equation*}
Recall that $\psi_{T}(\cdot) = \psi(\cdot/T)$ is a time cut-off function satisfying
\begin{equation*}
	\psi\in C_{0}^{\infty}(\mathbb R),\; \operatorname{supp}\psi \subseteq [\tau -1,\tau +1],\;
	\psi \equiv 1 \text{ on } \left[\tau - \frac{1}{2},\tau + \frac{1}{2}\right].
\end{equation*}
Then, we get that, for $T$ sufficiently small,
\begin{equation*}
	\|w\|_{X^{0,b}} \leq C\|w(\tau)\|,
	%\label{g5}
\end{equation*}
where $1/2<b<25/48$. By change the of initial time we can write
\begin{equation}
	\|w\|_{X^{0,b}_{[t_0-T,t_0+T]}} \leq C\|w(t_0)\|
	\label{g6}
\end{equation}
for every $t_0$. Now, taking the inner product of equation \eqref{g3}
with $2w$ in $L^2(\mathbb R)$, we find that
\begin{equation}
	\frac{d}{dt}\|w(t)\|^2 + 2\lambda\|w(t)\|^2 + 2\mu \|w_x(t)\|^2 =
	(u_x(t), w^2(t)).
	\label{g7}
\end{equation}
Multiplying \eqref{g7} by $e^{\lambda t}$ yields
\begin{equation}
	\frac{d}{dt}\left(e^{\lambda t}\|w(t)\|^2\right) + 2\mu  e^{\lambda
		t}\|w_x(t)\|^2 = e^{\lambda t}\left((u_x(t),w^2(t)) - \lambda
	\|w(t)\|^2\right).
	\label{g8}
\end{equation}
We integrate \eqref{g8} from $\tau$ to $t$ to obtain
\begin{equation}
	\begin{aligned}
		\|w(t)\|^2 &+ 2\mu \int_{\tau}^{t}e^{-\lambda(t-s)}\|w_x(s)\|^2\,ds\\
		&= \|w(\tau)\|^2e^{-\lambda(t-\tau)}+ \int_{\tau}^{t}e^{-\lambda(t-s)}
		\left((u_{x}(s), w^2(s)) - \lambda\|w(s)\|^2\right)ds.
	\end{aligned}
	\label{g9}
\end{equation}
We aim now to prove that the last term on the right hand side of \eqref{g9} is negative for large enough $N>0$, which implies the decay of $w$,
\begin{equation*}
	\|w(t)\|^2 \leq \|w(\tau)\|^2e^{-\lambda (t-\tau)}.
\end{equation*}
Consider the truncation function $\bar{\psi} = \psi_{[\tau,t]}$ and
the characteristic function $\chi = \chi_{[\tau,t]}$. Let $b$ and $b'$
such that $1/2<b<25/48$ and $1-b<b'<1/2$. Then, by duality, we find
that
\begin{equation}
	\begin{aligned}
		\int_{\tau}^{t}e^{-\lambda(t-s)}(u_x(s), w^2(s))\,ds
		&= \int_{-\infty}^{\infty}e^{-\lambda(t-s)}\chi(s)(\bar{\psi}(s)u(s), 
		\partial_x((\bar{\psi}(s)w(s))^2))\,ds\\
		&\leq \|\chi\bar{\psi}u\|_{X^{0,b'}}\|\partial_x(\bar{\psi}w)^2\|_{X^{0,-b'}}.
	\end{aligned}
	\label{g10}
\end{equation}
By the choice of $b, b'$ and by applying the bilinear estimate
\eqref{est2}, we have
\begin{equation}
	\|\partial_x(\bar{\psi}w)^2\|_{X^{0,b'}} \leq \|\partial_x(\bar{\psi}w)^2\|_{X^{0,b-1}}
	\leq \frac{C}{\sqrt{N}} \,\|w\|^2_{X^{0,b}_{[\tau,t]}}.
	\label{g11}
\end{equation}
On the other hand, since the characteristic function $\chi \in H_{t}^{b'}\cap L_{t}^{\infty}$ (recall $b' < 1/2$) and
$\bar{\psi}u$ belong to $X^{0,b}$ which is contained in
$L_{t}^{\infty}L_{x}^{2}$ (see \eqref{e5a}), one shows that
\begin{equation}
	\|\chi\bar{\psi}u\|_{X^{0,b'}} \leq
	\|\chi\|_{L^{\infty}_{t}}\|\bar{\psi}u\|_{X^{0,b'}} +
	\|\chi\|_{H_{t}^{b'}}\|\bar{\psi}u\|_{X^{0,b}} \leq
	C(\chi)\|u\|_{X^{0,b}_{[\tau,t]}},
	\label{g12}
\end{equation}
where we have used the definition of the Bourgain spaces. Here,
$C(\chi)$ depends on $\chi$, but is independent of $N$. Applying
\eqref{g11} and \eqref{g12} in \eqref{g10} yields
\begin{equation}
	\begin{aligned}
		\int_{\tau}^{t}e^{-\lambda(t-s)}(u_x(s), w^2(s))\,ds&\leq
		\frac{C(\chi)}{\sqrt{N}}\,\|u\|_{X^{0,b}_{[\tau,t]}}\|w\|^{2}_{X^{0,b}_{[\tau,t]}}\\ &\leq
		\frac{C(\chi)}{\sqrt{N}}\,\|u\|_{X^{0,b}_{[\tau,t]}}\frac{1}{t-\tau}
		\int_{\tau}^{t}\|w\|_{X^{0,b}_{[\tau,t]}}^{2}\,ds\\ &\leq
		\frac{C(\chi)}{\sqrt{N}}\,\|u\|_{X^{0,b}_{[\tau,t]}}\frac{1}{t-\tau}
		\int_{\tau}^{t}\|w(s)\|^2\,ds
	\end{aligned}
	\label{g13}
\end{equation}
here we have used \eqref{g6} at the last step. Thus, from \eqref{g9}
\begin{equation}
	\|w(t)\|^2 \leq \|w(\tau)\|^2e^{-\lambda(t-\tau)} +
	\int_{\tau}^{t}\left(\frac{C(\chi)}{\sqrt{N}(t-\tau)}\|u\|_{X^{0,b}_{[\tau,t]}}
	- \lambda e^{-\lambda(t-s)}\right)\|w(s)\|^2\,ds
	\label{g14}
\end{equation}
for $t$ is close enough to $\tau$. Naturally, this argument can be repeated to obtain \eqref{g14} for all $t\geq \tau$.
Hence, for $N$ large enough,
\begin{equation}
	\|w(t)\|^2 \leq \|w(\tau)\|^2e^{-\lambda(t-\tau)} = \|Q_Nu_{\tau}\|^2e^{-\lambda(t-\tau)}.
	\label{g15}
\end{equation}
Letting $\tau \rightarrow -\infty$, we deduce that $w$ decays
exponentially.
\subsubsection{Regularity of $v$}\label{regularity_v}
First, $v$ is defined globally in time, since $v = u - w$ and $u,\, w$
are defined globally in time. Next, we will prove the boundedness of
$v$ in $H^3(\mathbb R)$.

For $y(t) = P_Nv(t) = P_Nu(t)$, we have that $y$ is smooth and satisfies the estimate, thanks to Proposition \ref{absorbingset},
\begin{equation*}
	\lim_{\tau\rightarrow -\infty}\|y(t)\|_{H^{3}(\mathbb R)} \leq r(t)N^{3},
\end{equation*}
where $r(t)$ is the radius of the pullback absorbing set defined in \eqref{r}.
%Proposition~\ref{absorbingset}.

It remains to show that $Z = Q_Nv$, which solves
\begin{equation*}
	Z_t + Z_{xxx} + Q_N\left((y+Z)(y+Z)_{x}\right) + \lambda Z - \mu  Z_{xx} = Q_Ng(t)
	%\label{g16}
\end{equation*}
with $Z(\tau) = 0$, is bounded in $H^{3}(\mathbb R)$. 

It is easy to see that $Z = Q_Nv$
remains bounded in $L^2(\mathbb R)$. Now, by using the Agmon's inequality $\|f\|_{L^{\infty}(\mathbb R)} \leq \|u\|^{1/2}\|\partial_xu\|^{1/2}$
which implies 
\begin{equation*}
	\|f^2\|_{H^{-1}(\mathbb R)} \leq \|f\|^2
\end{equation*}
we get
\begin{equation}
	\begin{aligned}
		\|Z_{xxx}\| &\leq \|Z_t\| + \frac 12\|(y+Z)^2\|_{H^{-1}(\mathbb R)} +
		\lambda\|Z\| + \mu \|Z_{xx}\| + \|g(t)\|\\
		&\leq \|Z_t\| + \frac 14\|y\|^2 + \frac 14\|Z\|^2 +
		\lambda\|Z\| + \mu \|Z_{xx}\| + \|g(t)\|.
	\end{aligned}
	\label{g17}
\end{equation}

We observe from \eqref{g17} and interpolation inequalities that it is sufficient to show the boundedness of $Z_t$ in $L^2(\mathbb R)$. The equation for $Z'(t) = Z_t(t)$ is
\begin{equation}
	Z'_{t} + Z'_{xxx} + Q_{N}\partial_{x}((y+Z)Z') + \lambda Z' - \mu  Z'_{xx} = -Q_N\partial_x((y+Z)y_t) + g'(t)
	\label{g18}
\end{equation}
with initial condition $Z'(\tau) = Q_Ng(\tau) - Q_Ny(\tau)y_{x}(\tau)$
in $L^2(\mathbb R)$. Multiplying equation \eqref{g18} by $2Z'$,
integrating over $\mathbb R$, and using Young's inequality, we obtain
\begin{equation}\label{g19}
	\begin{aligned}
		\frac{d}{dt}\|Z'\|^2 &+ \frac{3\lambda}{2}\|Z'\|^2 + 2\mu \|Z'_{x}\|^2 \\
		&\leq C\left(\|\partial_x(y_t(y+Z))\|^2 + \int_{\mathbb R}(y+Z)\partial_{x}(Z')^2\,dx 
		+ \|g'(t)\|^2\right).
	\end{aligned}
\end{equation}
Using Young's inequality, then inverse (Poincar\'{e}-type)
inequalities (for functions with bounded spatial frequency) we find
that
\begin{equation*}
	\begin{aligned}
		\|\partial_{x}(y_{t}(y+Z))\|^2 \leq \|y_t\|^2_{H^1}\|y+Z\|^2_{H^1} \leq
		C(N)\left(1+ \|g(t)\|^2 + \|Z'\|^2\right),
		%\label{g20}
	\end{aligned}
\end{equation*}
where
\begin{equation*}
	\begin{aligned}
		\|\partial_{x}(y_{t}(y+Z))\|^2 &= \|\partial_{x}y_{t}(y+Z) + y_{t}\partial_{x}(y+Z)\|^2,\\
		&\leq \|\partial_{x}y_{t}\|\|y+Z\|_{L^\infty} + \|y_{t}\|_{L^\infty}\|\partial_{x}(y+Z)\|,\\
		&\leq 2\|y_{t}\|_{H^1}\|y+Z\|_{H^1}.
	\end{aligned}
\end{equation*}
Hence, \eqref{g19} implies that
\begin{equation*}
	\begin{aligned}
		\|Z'(t)\|^2 &\leq \|Z'(\tau)\|^2e^{-\lambda(t-\tau)} \\ &+
		C\int_{\tau}^{t}e^{-\lambda(t-s)}\int_{\mathbb
			R}\left((y+Z)\partial_x((Z')^2) + 1 + |g(s)|^2+
		|g'(s)|^2\right)\,dxds.
	\end{aligned}
	%\label{g21}
\end{equation*}
We proceed as in \eqref{g10}--\eqref{g15}, using the bilinear estimate
\eqref{est2} and the local well-posedness of \eqref{g18} in $L^2$, to
obtain
\begin{equation*}
	\|Z'(t)\|^2 \leq C\left(1 +
	\int_{-\infty}^{t}e^{-\lambda(t-s)}\left(\|g(s)\|^2+\|g'(s)\|^2\right)ds +
	\frac{1}{\sqrt{N}}\,\|Z'\|_{X^{0,b}_{[\tau,t]}}^2\right),
	%\label{g22}
\end{equation*}
for $t-\tau$ small enough, $1/2 < b < 25/48$, and consequently
\begin{equation}\label{twostar}
	\|Z'(t)\|^{2} \leq C\left(1+ \int_{-\infty}^{t}e^{-\lambda
		(t-s)}\left(\|g(s)\|^2+\|g'(s)\|^2\right)ds\right), \quad \forall
	\tau \leq \tau_0,
	%\label{g23}
\end{equation}
for some $\tau_0 \leq t$. This gives the boundedness of $v$ in $H^3(\mathbb R)$.

\begin{proposition}\label{boundedness}
	The pullback attractor $\mathcal A_{\mu } = \{A_\mu (t): t\in\mathbb
	R\}$ for\/ \textup{Eq.~\eqref{e1}} is bounded in $H^3(\mathbb R)$,
	i.e., for any $t\in\mathbb R$, $A_\mu (t)$ is a bounded subset of
	$H^3(\mathbb R)$.
\end{proposition}
\begin{proof}
	Let $\tau_n\rightarrow -\infty$ and let $\{u_{\tau_n}\}$ be a bounded
	sequence in $L^2(\mathbb R)$. Let $v_n$ and $w_n$ be the solutions of \eqref{g1} and \eqref{g2} respectively
	with initial data $u_{\tau_n}$. From \ref{regularity_v}, we have
	\begin{equation*}
		\{v_m(t, \tau_n+\cdot)\}_{n} \text{ is bounded in } C([-T,T];H^3(\mathbb R)),
		%\label{h1}
	\end{equation*}
	and, for the time derivative,
	\begin{equation*}
		\left\{ \frac{dv_{m}}{dt}(t,\tau_n+\cdot) \right\}_{n} \text{ is bounded
			in } C([-T,T];L^2(\mathbb R)),
		%\label{h2}
	\end{equation*}
	for each $T>0$ (and starting with $n$ sufficiently large so that
	$t-T\geq \tau_n$). Taking the inner product of \eqref{g1} with $2v$ in
	$L^2(\mathbb R)$ gives us an energy equality for $v$,
	\begin{multline*}
		\|v(t)\|^2 + 2\mu \int_{\tau}^{t}e^{-2\lambda(t-s)}\|v_x(s)\|^2\,ds
		\\ 
		= e^{-2\lambda(t-\tau)}\|v(\tau)\|^2+\, \int_{\tau}^{t}e^{-2\lambda(t-s)}\left(2(g(s),v(s)) + 
		(P_N(2v(s)w(s) + w^2(s)), v_x(s))\right)ds.
		%\label{g24}
	\end{multline*}
	Hence,
	\begin{multline}
		\|v(t)\|^2 + 2\mu \int_{t-T}^{t}e^{-2\lambda (t-s)}\|v_x(s)\|^2\,ds
		= e^{-2\lambda T}\|v(t-T)\|^2\\
		+ \, \int_{t-T}^{t}e^{-2\lambda (t-s)}\left(2(g(s),v(s)) 
		+ (P_N(2v(s)w(s) + w^2(s)), v_x(s))\right)ds.
		\label{h3}
	\end{multline}
	Note that $\{v_n\}$ is bounded in $H^3(\mathbb R)$ and $w_n
	\rightarrow 0$ in $L^2(\mathbb R)$ as $n\rightarrow 0$. Thus, by using
	\eqref{h3}, $v_n \rightarrow v_0$ in $L^2(\mathbb R)$ for some $v_0$
	in $H^3(\mathbb R)$ (since $\{v_n\}$ is bounded in $H^3(\mathbb
	R)$). On the other hand, by Theorem \ref{L2},
	\begin{equation*}
		u_n = v_n + w_n \rightarrow w_0 \enspace 
		\text{ in } L^2(\mathbb R) \text{ as } n\rightarrow \infty.
	\end{equation*}
	This implies that $v_0 = w_0$ and thus $w_0$ belongs to $H^3(\mathbb
	R)$. Hence we can conclude that the pullback attractor $\mathcal
	A_\mu $ is bounded in $H^3(\mathbb R)$.
\end{proof}

\subsection{Compactness in $H^{3}(\mathbb R)$}\label{regularity} 
Thanks to boundedness in $H^3(\mathbb R)$ and the compactness in $L^2(\mathbb R)$, we can use interpolation inequalities to get that the pullback attractor $\mathcal{A}_{\mu}$ for
\eqref{e1} is compact in $H^s(\mathbb R)$ for all $s<3$, and it is bounded
in $H^3(\mathbb R)$. We will prove that $\mathcal A_\mu $ is indeed
compact in $H^3(\mathbb R)$. To this end,  we
work with equations satisfied by $u_n' = du_n/dt$ and show, by
using the energy equation method applied to $u_n'$, some compactness
property of the $u_n'$ and then deduce the compactness of pullback
attractor in $H^3(\mathbb R)$.
\begin{theorem}\label{compactness}
The pullback attractor $\mathcal A_{\mu }$ corresponding to problem
\eqref{e1} is compact in $H^3(\mathbb R)$.
\end{theorem}
\begin{proof}
We can now restrict the process under study to $\mathcal A_\mu $. Since
$\mathcal A_\mu $ is bounded in $H^3(\mathbb R)$, we take initial data
$\{u_{\tau_n}\}$ being a bounded sequence in $H^3(\mathbb R)$, where
$\tau_n\rightarrow -\infty$. From equation \eqref{e1} and the fact
that $\mathcal A_{\mu }$ is compact in $H^s(\mathbb R)$ for $s<3$, it
is sufficient to prove the compactness of the sequence $\displaystyle
\left\{u'_n = \frac{du_n}{dt}\right\}$ in $L^2(\mathbb R)$.

We use the equation for $u'_n$,
\begin{equation}
\frac{du'_n}{dt} + \partial_{xxx}u'_n + (u_nu'_n)_x + \lambda u'_n
-\mu  \partial_{xx}u'_n = g'(t).
\label{h4}
\end{equation}
Since $\{u_n(t))\}_{n}$ is bounded in $H^3(\mathbb R)$ and converges
strongly to $u_0(t)$ in $C([\tau-T,\tau+T];L^2(\mathbb R))$, we see
that
\begin{equation}
\begin{aligned}
u'_n(t) \rightarrow u'_0(t) &\text{ strongly in } C([-T,T];H^{-s}(\mathbb R))\\
&\text{ and weakly star in } L^{\infty}([-T,T];L^2(\mathbb R)).
\end{aligned}
\label{h5}
\end{equation}
Now, consider the energy equality for $u'_n$ by multiplying \eqref{h4} with $u'_n$ in $L^2(\mathbb R)$ then integrating over $(t-T, t)$
\begin{multline}
\|u'_n(t)\|^2 + \mu \int_{t-T}^{t}e^{-2\lambda(t-s)}\|\partial_{x}u'_n(s)\|^2\,ds
= e^{-2\lambda T}\|u'_n(t-T)\|^2\\
-\,\int_{t-T}^{t}e^{-2\lambda(t-s)}(\partial_xu_n, u_n'^2)\,ds
+ \int_{t-T}^{t}e^{-2\lambda(t-s)}(g'(s),u'_n(s))\,ds.
\label{h6}
\end{multline}
We want to pass to the limit as $n\rightarrow +\infty$ in \eqref{h6}. On one hand, using
\eqref{h5}, as $n\rightarrow \infty$, we get
\begin{equation*}
\int_{t-T}^{t}e^{-2\lambda(t-s)}\left(\partial_{x} u_n - \partial_xu_0, u_0'^2\right)ds 
\rightarrow 0.
%\label{h7}
\end{equation*}
On the other hand, for $0<s\leq 23/48$ and $1/2<b<25/48$, we take
$b''$ such that $1-b<b''<1/2$ and proceed as in
\eqref{g10}--\eqref{g13} (but with $-s$ and using Lemma
\ref{bilinear2}) to find
\begin{equation*}
\left|\int_{t-T}^{t}e^{-2\lambda(t-s)}\left(\partial_xu_n, u_n'^2 - u_0^2\right)ds\right| 
\leq c\|u_n\|_{X^{s,b}_{[\tau,\tau+T]}}\|u_n'-u_0'\|^2_{X^{-s,b}_{[\tau,\tau+T]}}.
%\label{h8}
\end{equation*}
From \eqref{h5} and the well-posedness of \eqref{h4} in
$X_{[\tau,\tau+T]}^{-s,b}$, we can pass to the limit as $n$ goes to
infinity to obtain
\begin{equation*}
\lim_{n\rightarrow \infty}\int_{t-T}^{t}e^{-2\lambda(t-s)}
\left(\partial_xu_n(t+s),u_n'(t+s)^2\right)ds
=\int_{t-T}^{t}e^{-2\lambda(t-s)}\left(\partial_xu_0(s), u_0'(s)^2\right)ds.
%\label{h9}
\end{equation*}
By \eqref{h5} once again, we get
\begin{equation*}
\lim_{n\rightarrow \infty}\int_{t-T}^{t}e^{-2\lambda(t-s)}
\left(g'(s),u_n'(s)\right)ds 
= \int_{t-T}^{t}e^{-2\lambda(t-s)}\left(g'(s),u_0'(s)\right)ds.
%\label{h9_1}
\end{equation*}
Then, from \eqref{h6},
\begin{multline*}
\limsup_{n\rightarrow \infty}\left(\|u_n'(t)\|^2
+ \mu \int_{t-T}^{t}e^{-2\lambda(t-s)}\|\partial_{xx}u'_n(s)\|^2ds\right)\\
\leq Ce^{-2\lambda T}-\int_{t-T}^{{t}}e^{-2\lambda(t-s)}
\left(\partial_xu_0(s), u_0'(s)^2\right)ds
+ \int_{t-T}^{t}e^{-2\lambda(t-s)}\left(g'(s),u_0'(s)\right)ds.
%\label{h10}
\end{multline*}
By employing a corresponding $L^2$ energy equality for $u_0'$
and using
\begin{equation*}
\limsup_{n\rightarrow \infty}\mu \int_{t-T}^{t}e^{-2\lambda(t-s)}
\|\partial_{x}u'_n(s)\|^2\,ds\leq \int_{t-T}^{t}e^{-2\lambda(t-s)}
\|\partial_{x}u_0'(s)\|^2\,ds,
%\label{h11}
\end{equation*}
we find that
\begin{equation*}
\limsup_{n\rightarrow\infty}\|u_n'(t)\|^2 \leq Ce^{-2\lambda T} + \|u_0'(t)\|^2.
%\label{h12}
\end{equation*}
We let $T$ go off to infinity to obtain
\begin{equation*}
\limsup_{n\rightarrow \infty}\|u'_n(t)\|^2 \leq \|u_0'(t)\|^2.
%\label{h13}
\end{equation*}
Therefore,
\begin{equation*}
u'_n(t) \rightarrow u_0'(t) \enspace \text{ strongly in } L^2(\mathbb R).
%\label{h14} 
\end{equation*}
This complete the proof.
\end{proof}

%\section{Pullback attractors for KdV equation}
%In this section, we show that there exists a pullback attractor $\mathcal A = \{A(t)\}$ for the KdV equation, that is the equation \eqref{e1} with $\mu=0$
%
%and moreover, this attractor is compact in $H^3(\mathbb R)$.
%\begin{theorem}
%	Assume that the external force $g$ satisfies \eqref{e2}, then there exists a unique pullback attractor $\mathcal A = \{A(t)\}$ in $L^2(\mathbb R)$ for the equation \eqref{KdV}. Moreover, $A(t)$ is compact in $H^3(\mathbb R)$ for all $t\in \mathbb R$.
%\end{theorem}
%\begin{proof}
%\end{proof}
\section{Upper semicontinuity at $\mu  = 0$}\label{sec:5}
In this section, we will prove that, as $\mu  \rightarrow 0^+$, the
pullback attractor $\mathcal A_{\mu }$ of the KdV-Burgers equation
converges to the pullback attractor $\mathcal A_0$ of the weakly
damped KdV equation in the sense of the Hausdorff semi-distance. To this end, we first show
\begin{theorem}\label{thrm:KdV}
	Assume that the external force $g$ satisfies \eqref{e2}, then there exists a unique pullback attractor $\mathcal A_0 = \{A_0(t)\}$ in $L^2(\mathbb R)$ for the equation
	\begin{equation}\label{KdV}
		\begin{cases}
			u_t + u_{xxx} + uu_x + \lambda u = g(t,x), &x\in\mathbb R, \quad t>0,\\
			u(\tau,x) = u_\tau(x), &x\in\mathbb R.
		\end{cases}
	\end{equation}
	Moreover, $A_0(t)$ is compact in $H^3(\mathbb R)$ for all $t\in \mathbb R$.
\end{theorem}
\begin{proof}
	The arguments of the proof of this Theorem follow from the previous section, except the fact that we need to work in the Bourgain spaces corresponding to the free equation of \eqref{KdV} instead of the spaces defined by \eqref{BourgainSpace}. More precisely, we use the Bourgain space
	\begin{equation*}
		\mathcal X^{s,b} = \{u: \mathbb R\rightarrow \mathbb R: \|u\|_{\mathcal X^{s,b}} < +\infty \}
	\end{equation*}
	where
	\begin{equation*}
		\|u\|_{\mathcal X^{s,b}} = \left\|\langle \tau - \xi^3\rangle\langle \xi\rangle^s \widehat{u}(\tau,\xi)\right\|_{L^2(\mathbb R^2)}.
	\end{equation*}
	The crucial estimate \eqref{est1} in this case is replaced by
	\begin{equation*}
		\|D(u^2)\|_{\mathcal X^{s,b'-1}} \leq c_{s}\|u\|_{\mathcal X^{s,b}}\|u\|_{\mathcal{X}^{-1/2,b}},
	\end{equation*}
	which was proved in \cite[Proposition 2.2]{Goubet}.
	
	Repeat the arguments in the previous section we first get the existence of the pullback attractor $\mathcal A_0 = \{A_0(t)\}$ in $L^2(\mathbb R)$, then by splitting solution technique (see \eqref{g1} and \eqref{g2}) we obtain the compactness of $\mathcal A_0$ in $H^3(\mathbb R)$.
\end{proof}
Thanks to the compactness of $\mathcal A_{\mu}$ and $\mathcal A_0$ in $H^3(\mathbb R)$, we can from now on assume that the initial data $u_\tau$
belongs to $H^3(\mathbb R)$. We also additionally assume that the external force $g$
belongs to $L^2_{\text{loc}}(\mathbb R;H^1(\mathbb R))$ and satisfies
\begin{equation}
\lim_{t\rightarrow -\infty} e^{-\lambda t}\int_{-\infty}^{t}e^{\lambda
  s}\|g(s)\|_{H^1(\mathbb R)}^2\,ds < +\infty.
\label{Gbis}
\end{equation}

We need the following lemmas.

\begin{lemma}[\cite{Robinson}]\label{lem-2}
Let $V$ be a Banach space, $H$ be a Hilbert space, and $V\subseteq H$
be continuous. Suppose that $\{u_n\}$ is bounded in
$L^{\infty}(\tau,T;V)$,
\begin{equation*}
\underset{t\in[\tau,T]}{\operatorname{ess}\sup}\,\|u_n(t)\|_{V}\leq C,
%\label{003}
\end{equation*}
and that $u_n\rightharpoonup u$ in $L^2(\tau,T;V)$, where, in addition,
$u\in C([\tau,T];H)$. Then
\begin{equation*}
\sup_{t\in[\tau,T]}\|u(t)\|_{V}\leq C.
%\label{004}
\end{equation*}
\end{lemma}
\begin{lemma}\label{lem-1}
Let $u(t)$ be a solution of \eqref{e1} with initial data $u_{\tau}\in
H^3(\mathbb R)$. Then $u$ is bounded in $L^{\infty}(\tau,t;H^2(\mathbb
R))$ uniformly in $\mu \in[0,1]$, i.e.,
\begin{equation}
\|u\|_{L^{\infty}(\tau,t;H^2(\mathbb R))} \leq C,
\label{000}
\end{equation}
where $C$ depends only on $\tau,\,t$, and $\|u_{\tau}\|_{H^{3}(\mathbb
  R)}$, but is independent of $\mu $.
\end{lemma}
\begin{proof}
As in \cite{DlotkoSun}, we consider the regularized equation
\begin{equation*}
\begin{cases}
\enspace u_t^{\epsilon} + \epsilon u^{\epsilon}_{xxxx} + u^{\epsilon}_{xxx} +
u^{\epsilon}u^{\epsilon}_{x} + \lambda u^{\epsilon} - \mu 
u^{\epsilon}_{xx} = g(t),\\ 
\enspace u|_{t=\tau} = u_{\tau}\in H^3(\mathbb R).
\end{cases}
%\label{005}
\end{equation*}
By results of \cite{DlotkoSun}, 
\begin{equation*}
\{u^{\epsilon}\}_{\epsilon\in (0,1]} \text{ is bounded in } 
L^{\infty}(\tau,t;H^2(\mathbb R)) \text{ uniformly in } \mu 
%\label{006}
\end{equation*}
and
\begin{equation*}
u^{\epsilon}\rightharpoonup u \text{ in } 
L^{2}(\tau,t;H^2(\mathbb R)) \text{ as } \epsilon \rightarrow 0^+.
%\label{007}
\end{equation*}
Using Lemma \ref{lem-2}, we immediately get \eqref{000}.
\end{proof}

\begin{lemma}\label{lem0}
For any $t\in \mathbb R$, the set 
\begin{equation*}
\bigcup_{\tau\leq t}\bigcup_{\mu \in[0,1]} A_{\mu }(t)
\end{equation*}
is bounded in $L^2(\mathbb R)$.
\end{lemma}
\begin{proof}
From Proposition~\ref{absorbingset}, $\{U_\mu (t,\tau)\}$ has a
pullback absorbing set $\mathcal B_\mu  =$ \linebreak
$\{B_\mu (t)\}_{t\in\mathbb R}$ in $L^2(\mathbb R)$, where $B_\mu (t)$
is the ball centered at the origin with radius
\begin{equation*}
R(t) = \left(\frac{e^{-\lambda
    t}}{\lambda}\int_{-\infty}^{t}e^{\lambda
  s}\|g(s)\|^2\,ds\right)^{1/2}.
\end{equation*}
Note that $R(t)$ is independent of $\mu $. Combining this with
$A_\mu (t)\subseteq B_{\mu }(t)$, we get
\begin{equation}\label{ll1}
\bigcup_{\mu \in[0,1]}A_\mu (t)\subseteq B_0(t).
\end{equation}
On the other hand, by assumption \eqref{Gbis} on $g$, 
\begin{equation*}
\limsup_{t\rightarrow -\infty}R(t) < +\infty.
\end{equation*}
This implies that
\begin{equation}\label{ll2}
\bigcup_{\tau\leq t}B_0(\tau) \text{ is bounded in } L^2(\mathbb R).
\end{equation}
From \eqref{ll1} and \eqref{ll2}, we get the desired result.
\end{proof}

\begin{lemma}\label{lem1}
For any $t\in \mathbb R$, the set $\displaystyle 
\bigcup\limits_{\mu \in[0,1]}A_{\mu }(t)$ is precompact in $L^2(\mathbb R)$.
\end{lemma}
\begin{proof}
By the invariance property of $\mathcal A_{\mu }$, 
\begin{equation*}
A_{\mu }(t) = U_{\mu }(t,\tau)A_{\mu }(\tau),\quad \forall t\geq \tau.
\end{equation*}
In the previous section, we proved that $A_\mu (t)$ is compact in
$H^3(\mathbb R)$. Thus, any element $u(t)$ of $\mathcal A_\mu (t)$ can
be considered as a solution of \eqref{e1} at time $t$ with initial
data belonging to $H^3(\mathbb R)$. Hence, the result follows directly
from \cite[Theorem 5.2]{AB12}. For the reader's convenience, we
sketch the proof here. Namely, we prove that, for any $\eta>0$, there
exists a finite covering of the set $\bigcup_{\mu \in(0,1]}A_{\mu }(t)$
  in $L^2(\mathbb R)$ by balls of radius less than $\eta$. First,
  using the boundedness of $u$ from Lemma \ref{lem-1}, we multiply
  equation \eqref{e1} by $\theta\left(x^2/k^2\right)u$, where
  $\theta\colon [0,\infty) \rightarrow [0,1]$ is the smooth cut-off
  function satisfying $\theta|_{[0,1]} = 0$ and $\theta|_{[2,+\infty)}
    = 1$, to obtain that
\begin{equation}
\int_{|x|\geq k}|u|^2\,dx \leq \frac{\eta}{4},\quad \forall u\in 
\bigcup_{\mu \in (0,1]}A_{\mu }(t),
\label{008}
\end{equation}
for some $k>0$.

Next, from the boundedness of solutions in $H^2(\mathbb R)$, we conclude that
\begin{equation}
\|u\|_{H^1(-k,k)} \leq C, \quad\forall u\in
\bigcup_{\mu \in(0,1]}A_{\mu }(t).
\label{009}
\end{equation}
Because $H^1(-k,k)\hookrightarrow L^2(-k,k)$ compactly, we obtain from
\eqref{008} and \eqref{009} that 
$\bigcup_{\mu \in(0,1]}A_{\mu }(t)$ has the desired finite covering in
$L^2(\mathbb R)$ by balls of radius less than~$\eta$. This completes
the proof.
\end{proof}

The following result is the key in proving the upper semicontinuity.
\begin{lemma}\label{lem2}
Let $\{x_n\}\subset H^3(\mathbb R)$ and $\{\mu _n\}\subset [0,1]$ be
two sequences such that $\{x_n\}$ is bounded in $H^3(\mathbb R)$, $x_n
\rightarrow x_0$ in $L^2(\mathbb R)$ and $\mu _n\rightarrow 0$. Then,
for any $t\geq \tau$, 
\begin{equation}
\|U_{\mu _n}(t,\tau)x_n - U_{0}(t,\tau)x_0\| \rightarrow 0 \enspace \textup{ as
} n\rightarrow 0.
\label{j1}
\end{equation}
\end{lemma}
\begin{proof}
Denote $u_n(t) = U_{\mu _n}(t,\tau)x_n$ and $v(t) =
U_0(t,\tau)x_0$. Then, $u_n$ and $v$ solve the equations
\begin{equation*}
\begin{cases}
\enspace \partial_tu_n + \partial_{xxx}u_n + u_n\partial_xu_n +
\lambda u_n - \mu _n\partial_{xx}u_n = g(t,x)\\ \enspace u_n(\tau) = x_n,
\end{cases}
%\label{j2}
\end{equation*}
and
\begin{equation*} %\label{j3}
\begin{cases}
\enspace v_t + v_{xxx} + vv_x + \lambda v = g(t,x)\\
\enspace v(\tau) = x_0,
\end{cases}
\end{equation*}
respectively. Let $w_n(t) = u_n(t) - v(t)$. Then we have
\begin{equation}
\partial_tw_n + \partial_{xxx}w_n + (u_n\partial_xu_n-vv_x) + \lambda w_n - \mu _n\partial_{xx}u_n = 0
\label{j4}
\end{equation}
with initial data $w_n(\tau) = x_n - x_0$. Multiplying \eqref{j4} by
$w_n$ and integrating over $\mathbb R$, we obtain
\begin{equation}
\frac{1}{2}\frac{d}{dt}\|w_n\|^2 + \int_{\mathbb
  R}(u_n\partial_xu_n-vv_x)w_n\,dx + \lambda\|w_n\|^2
+\mu _n\int_{\mathbb R}\partial_xu_n\partial_xw_n\,dx = 0.
\label{j5}
\end{equation}
Integrating by parts and using interpolation inequalities, we find that
\begin{equation}
\begin{aligned}
\left|\int_{\mathbb R}(u_n\partial_xu_n-vv_x)w_n\,dx\right| 
&= \frac{1}{2}\left|\int_{\mathbb R}(u_n^2 - v^2)_xw_n\,dx\right|\\
&=\frac{1}{2}\left|\int_{\mathbb R}(u_n+v)w_n\partial_xw_n\,dx\right|\\
&\leq\int_{\mathbb R}\left(|\partial_xu_n| + |v_x|\right)w_n^2\,dx\\
&\leq \left(\|\partial_xu_n\|_{L^{\infty}(\mathbb R)}+\|v_{x}\|_{L^{\infty}(\mathbb R)}\right)
\|w_n\|^2\\
&\leq \left(\|u_n\|_{H^2(\mathbb R)}+\|v\|_{H^2(\mathbb R)}\right)\|w_n\|^2.
\end{aligned}
\label{j6}
\end{equation}
On the other hand, by Young's inequality, 
\begin{equation}
\begin{aligned}
\left|\mu _n\int_{\mathbb R}\partial_xu_n\partial_xw_n\,dx\right|&\leq 
\mu _n\left|\int_{\mathbb R}\partial_xu_n\partial_xw_n\,dx\right|\\
&\leq \mu _n\|\partial_xu_n\|_{L^{2}(\mathbb R)}\|\partial_nw_n\|_{L^{2}(\mathbb R)}\\
&\leq \mu _n\|u_n\|_{H^2(\mathbb R)}\|w_n\|_{H^2(\mathbb R)}.
\end{aligned}
\label{j7}
\end{equation}
From \eqref{j5}--\eqref{j7}, we obtain
\begin{equation*}
\frac{d}{dt}\|w_n\|^2 \leq \left(\|u_n\|_{H^2(\mathbb R)}+\|v\|_{H^2(\mathbb R)}\right)
\|w_n\|^2 + \mu _n\|u_n\|_{H^2(\mathbb R)}\|w_n\|_{H^2(\mathbb R)}.
%\label{j8}
\end{equation*}
Integrating from $\tau$ to $t$ and using Lemma \ref{lem-1} yields
\begin{equation*}
\begin{aligned}
\|w_n(t)\|^2 &\leq \|w_n(\tau)\|^2 + \int_{\tau}^{t}
\left(\|u_n(s)\|_{H^2(\mathbb R)}+\|v(s)\|_{H^2(\mathbb R)}\right)\|w_n(s)\|^2\,ds\\
&\qquad+\mu _n
\int_{\tau}^{t}\|u_n(s)\|_{H^2(\mathbb R)}\|w_n(s)\|_{H^2(\mathbb R)}\,ds\\
&\leq \|w_n(\tau)\|^2 + 
\left(\|u_n\|_{L^{\infty}(\tau,t;H^2(\mathbb R))}+\|w_n\|_{L^{\infty}(\tau,t;H^2(\mathbb R))}
\right)\int_{\tau}^{t}\|w_n(s)\|^2\,ds\\
&\qquad +\mu _n\|u_n\|_{L^{2}(\tau,t;H^2(\mathbb R))}\|w_n\|_{L^{2}(\tau,t;H^2(\mathbb R))}.
\end{aligned}
%\label{j9}
\end{equation*}
Applying Gronwall's lemma, we find that
\begin{equation*}
\begin{aligned}
\|w_n(t)\|^2 &\leq \left(\|w_n(\tau)\|^2 + \mu _n\|u_n\|_{L^{2}(\tau,t;H^2(\mathbb R))}
\|w_n\|_{L^{2}(\tau,t;H^2(\mathbb R))}\right)\\
&\quad \times \exp\left(\left(\|u_n\|_{L^{\infty}(\tau,t;H^2(\mathbb R))}
+\|w_n\|_{L^{\infty}(\tau,t;H^2(\mathbb R))}\right)(t-\tau)\right).
\end{aligned}
%\label{j10}
\end{equation*}
Now let $n\rightarrow \infty$ and note that $w_n(\tau) = x_n - x_0$, then
use \eqref{000} and $x_n\rightarrow x_0$ in $L^2(\mathbb R)$, to get
the desired result.
\end{proof}

We can now prove the main result of this section.
\begin{theorem}\label{upper}
Assume that $g\in L^2_{\textup{loc}}\left(\mathbb R;H^1(\mathbb
R)\right)$ satisfies \eqref{Gbis}. Then the family of pullback
attractors $\{\mathcal A_\mu \}_{\mu \in[0,1]}$ corresponding to\/
\textup{Eq.~\eqref{e1}} is upper semi-continuous in $L^2(\mathbb R)$
at $\mu  = 0$, i.e.,
\begin{equation*}
\lim_{\mu \rightarrow 0^+}\operatorname{dist}_{L^2(\mathbb R)}(A_\mu (t), A_0(t)) = 0 
\quad \textup{ for all } t\in \mathbb R,
\end{equation*}
where $\operatorname{dist}_{L^2(\mathbb R)}$ denotes the Hausdorff
semi-distance in $L^2(\mathbb R)$.
\end{theorem}
\begin{proof}
We prove this by contradiction. Assume that there exists a $t_0\in
\mathbb R$ such that
\begin{equation*}
\operatorname{dist}_{L^2(\mathbb R)}(A_\mu (t_0), A_0(t_0)) \not\rightarrow 0 
\text{ as } \mu  \rightarrow 0^+.
\end{equation*}
By the compactness of $\mathcal A_{\mu }$, we can choose a positive
constant $\delta>0$, a sequence $\mu _n\rightarrow 0$, and a sequence
$\{\varphi_n\}$ such that $\varphi_n\in A_{\mu _n}(t_0)$ and
\begin{equation}
\operatorname{dist}_{L^2(\mathbb R)}(\varphi_n, A_0(t_0)) \geq \delta
\enspace \text{ for all } n\geq 1.
\label{j11}
\end{equation}
By Lemma \ref{lem0}, the set $\bigcup_{\tau\leq
  t_0}\bigcup_{\mu \in[0,1]}A_\mu (\tau)$ is bounded in $L^2(\mathbb
R)$. Then, by the pullback attracting property of $\{U_{0}(t,\tau)\}$,
we can choose a $\tau_0\leq t_0$ satisfying
\begin{equation*}
\operatorname{dist}_{L^2(\mathbb R)}(U_0(t_0,\tau_0)\psi, A_0(t)) \leq
\frac{\delta}{4} \enspace \text{ for all } \psi\in \bigcup_{\tau\leq
  t_0}\bigcup_{\mu \in[0,1]}A_\mu (\tau_0).
\end{equation*}
Since $\mathcal A_{\mu _n}$ is invariant, for each $n\geq 1$, 
\begin{equation*}
U_{\mu _n}(t_0,\tau_0)A_{\mu _n}(\tau_0) = A_{\mu _n}(t_0).
\end{equation*}
Thus, there exist $\psi_n\in A_{\mu _n}(\tau_0)$ such that
\begin{equation}
\varphi_n = U_{\mu _n}(t_0,\tau_0)\psi_n.
\label{j12}
\end{equation}
The sequence $\{\psi_n\}$ is contained in the set $\bigcup_{n\geq
  1}A_{\mu _n}(\tau_0)\subseteq
\bigcup_{\mu \in[0,1]}A_{\mu }(\tau_0)$. By Lemma \ref{lem1}, we can
assume that $\psi_n \rightarrow \psi_0$ in $L^2(\mathbb R)$. Applying
Lemma \ref{lem2}, we have
\begin{equation*}
\varphi_n = U_{\mu _n}(t_0,\tau_0)\psi_n \rightarrow
U_{0}(t_0,\tau_0)\psi_0 \enspace \text{ in } L^{2}(\mathbb R) \text{ as }
n\rightarrow \infty.
%\label{j13}
\end{equation*}
Choose $N$ so large that
\begin{equation}
\|U_{\mu _n}(t_0,\tau_0)\psi_n - U_{0}(t_0,\tau_0)\psi_0\| \leq \frac{\delta}{4} 
\enspace \text{ for all } n\geq N.
\label{j14}
\end{equation}
Then, by \eqref{j11}, \eqref{j12}, and \eqref{j14}, 
\begin{equation*}
\begin{aligned}
0< \delta &\leq \operatorname{dist}_{L^2(\mathbb
  R)}(U_{\mu _N}(t_0,\tau_0)\psi_N, A_0(t_0))\\
&\leq \|U_{\mu _N}(t_0,\tau_0)\psi_N - U_0(t_0,\tau_0)\psi_0\| +
\operatorname{dist}_{L^2(\mathbb R)}(U_0(t_0,\tau_0)\psi_0, A_0(t_0))\\ 
&\leq \frac{\delta}{4} + \frac{\delta}{4}= \frac{\delta}{2}.
\end{aligned}
%\label{j15}
\end{equation*}
This is a contradiction which thus completes the proof.
\end{proof}

\section{Finite fractal dimensions}
In this section, we consider the case $\mu>0$ and prove that  the pullback attractor $\{\mathcal{A}_{\mu }(t)\}_{t\in\mathbb R}$ has finite fractal dimensions in $H^1(\R)$. This will be achieved by proving and applying an extension of the splitting method, see e.g. \cite[Theorem 4.15]{CL08}. Therein the abstract results were written for global attractors of a semigroup. Here we prove an generalized version of that result for pullback attractors of non-autonomous processes. We start with the definition of fractal dimension: A subset $M$ of a Banach space $X$ is called to have finite fractal dimension in $X$ if 
\begin{equation*}
    \dim_f^{X}(M):= \limsup_{\varepsilon\to 0} \frac{\ln n(M,\varepsilon)}{\ln(1/\varepsilon)} < +\infty
\end{equation*}
where $n(M,\varepsilon)$ is the minimal number of closed balls in $X$ of radius $\varepsilon$ which covers the set $M$.
\begin{theorem}\label{new_prop}
    Let $X$ be a Banach space and $\{\mathcal{A}(t)\}_{t\in\R}$ be a family of compact sets in $X$. Assume that there exists a family of operators $\mathscr{V}(t, s): \mathcal{A}(s) \to X$, $t\ge s$, such that
    \begin{itemize}
        \item[(i)] $\mathcal{A}(t) \subset \mathscr{V}(t, s)\mathcal{A}(s)$, for $t \geq s$.
        \item[(ii)] $\mathscr{V}(t, s)$ is Lipschitz on $\mathcal{A}(s)$, i.e, there exists $L_{t, s} > 0$ such that
        \begin{equation*}
            \|\mathscr{V}(t, s)v_1-\mathscr{V}(t, s)v_2\| \leq L_{t, s}\|v_1-v_2\|,\quad \forall v_1, v_2 \in \mathcal A(s),
        \end{equation*}
        where the family of constants $\{L_{t,s}\}_{t\ge s}$ satisfies $$ \limsup_{s \to -\infty} L_{t, s} < +\infty.$$
        \item[(iii)] There exist compact semi-norms $n_1, n_2$ on $X$ such that
        \begin{align*}
            &\|\mathscr{V}(t, s)v_1-\mathscr{V}(t, s)v_2\|\\
            &\leq \eta_{t, s}\|v_1-v_2\| + K_{t, s}\Big(n_1(v_1-v_2)+n_2(\mathscr{V}(t, s)v_1-\mathscr{V}(t, s)v_2)\Big),
        \end{align*}
        for any $v_1, v_2 \in \mathcal{A}(s)$, where the family of constants $\eta{_{t,s}} > 0$, $K_{t, s} >0$ satisfy 
        $$ \limsup_{s \to -\infty} \eta_{t, s} \in (0,1/2) \quad \text{ and } \quad \limsup_{s \to -\infty}K_{t, s} < +\infty.$$
        \end{itemize}
        Then, $\{\mathcal{A}(t)\}_{t\in\R}$ has finite fractal dimensions in $X$, that is
        \begin{equation*}
            \dim_{f}^X(\mathcal{A}(t))  < +\infty \text{ for each } t\in \mathbb R.
        \end{equation*}
\end{theorem}
\begin{proof}
    Fix $t\in\mathbb{R}$. Define
    \begin{equation*}
    \tilde\eta_t:=\limsup_{s\to-\infty}\eta_{t,s},\qquad 
    \tilde K_t:=\limsup_{s\to-\infty}K_{t,s},\qquad 
    \tilde L_t:=\limsup_{s\to-\infty}L_{t,s}.
    \end{equation*}
    By assumptions $0<\tilde\eta_t<1/2$, $\tilde K_t<\infty$ and $\tilde L_t<\infty$. Choose  $\eta, K, L$ such that $\tilde{\eta}_t < \eta < 1/2$, $\tilde{K}_t<K$, and $\tilde{L}_t<L$. Then, there exists $S \leq t$ such that for all $s \leq S$,
    \begin{equation*}
        \eta_{t,s} \leq \eta, \qquad K_{t, s} \leq K, \qquad L_{t, s} \leq L,
    \end{equation*}
    and for every $v_1,v_2\in \mathcal{A}(s)$,
    \begin{equation}
    \|\mathscr{V}_{t,s}v_1-\mathscr{V}_{t,s}v_2\|
    \le \eta\|v_1-v_2\|+K\left[n_1(v_1-v_2)+n_2\left(\mathscr{V}_{t,s}v_1-\mathscr{V}_{t,s}v_2\right)\right].
    \label{iii}
    \end{equation}
    For any $s\leq S$ we define a pseudo-metric $\varrho_{t,s}$ on $\mathcal{A}(s)$ by
    \begin{equation*}
    \varrho_{t,s}(x,y):=K\Bigl[n_1(x-y)+n_2\bigl(\mathscr{V}_{t,s}x-\mathscr{V}_{t,s}y\bigr)\Bigr],\qquad x,y\in \mathcal{A}(s).
    \end{equation*}
    Because $n_1,n_2$ are compact semi-norms and $\mathcal{A}(s)$ is bounded, $\varrho_{t,s}$ is compact. From $(iii)$, we have
    \begin{equation}
        \|\mathscr{V}(t, s)x-\mathscr{V}(t, s)y\| \leq \eta \|x-y\| + \varrho_{t, s}(x, y), \quad \forall x, y \in \mathcal{A}(s).
    \label{lipschitz}
    \end{equation}
    Fix $\delta \in (0, 1-2\eta)$ and set $q:=2\eta+\delta<1$. For any $\varepsilon>0$ and any $B \subset \mathcal{A}(s)$ with $\diam B \leq 2\varepsilon$, let $m_{\varrho_{t, s}}(B, \delta\varepsilon)$ be the maximal cardinality of a $(\delta\varepsilon, \varrho_{t, s})$-distinguishable subset of $B$ (where a subset $Q \subseteq B$ is defined as $(\delta\varepsilon, \varrho_{t, s}$)-distinguishable if for any two distinct elements $x,y \in Q$, the pseudo-metric satisfies $\varrho_{t, s}(x,y)>\delta\varepsilon$). Define the local capacity
    \begin{equation*}
        \mathcal{C}_{\varrho_{t,s}}(\mathcal{A}(s); \varepsilon, \delta\varepsilon) := \sup\{\ln m_{\varrho_{t,s}} (B, \delta\varepsilon): B \subseteq \mathcal{A}(s), \diam(B) \leq 2\varepsilon\},
    \end{equation*}
    Now, we estimate $\mathcal{C}_{\varrho_{t,s}}$ . For any $x, y \in B$, using the Lipschitz property $(ii)$ with constant $L$,
    \begin{equation*}
        \|x-y\|^2 + \|\mathscr{V}_{t, s}x-\mathscr{V}_{t, s}y\|^2 \leq (1+L^2)\|x-y\|^2 \leq (1+L^2)(2\varepsilon)^2 = 4(1+L^2)\varepsilon^2.
    \end{equation*}
    Thus, the points $Z=(x, \mathscr{V}_{t, s}x)$ lie in a ball of radius $R=2\sqrt{1+L^2}\varepsilon$ in $X \times X$. The condition $\varrho_{t, s}(x, y) > \delta\varepsilon$  is equivalent to
    \begin{equation*}
        n_1(x-y)+n_2\left(\mathscr{V}_{t, s}x-\mathscr{V}_{t, s}y\right) > \frac{\delta\varepsilon}{K}.
    \end{equation*}
    By defining $\lambda=K/(\delta\varepsilon)$, the scaled points $\widetilde{Z}=\lambda Z$ lie in the ball of radius $\frac{K}{\delta\varepsilon}R = \frac{2K\sqrt{1+L^2}}{\delta}$ and satisfy $\eta_1(\widetilde{x}-\widetilde{y})+\eta_2(\mathscr{V}_{t,s}\widetilde{x}-\mathscr{V}_{t,s}\widetilde{y})>1$ for all $\widetilde{x}, \widetilde{y} \in \widetilde{Z}$. Hence the maximal number of such points is at most $m_0\!\left(\frac{2K\sqrt{1+L^2}}{\delta}\right)$, where $m_0(R)$ is the maximal cardinality of pairs $(x_i,y_i)\in X\times X$ with $\|x_i\|^2+\|y_i\|^2\le R^2$ and $\eta_1(x_i-x_j)+\eta_2(y_i-y_j)>1$ for $i\neq j$. Therefore
    \begin{equation*}
        m_{\varrho_{t,s}}(B,\delta\varepsilon)\le m_0\!\left(\frac{2K\sqrt{1+L^2}}{\delta}\right),
    \end{equation*}
    and taking the supremum over all such $B$ gives
    \begin{equation}
        \mathcal{C}_{\varrho_{t,s}}(\mathcal{A}(s);\varepsilon,\delta\varepsilon)
        \leq \ln m_0\!\left(\frac{2K\sqrt{1+L^2}}{\delta}\right).
    \label{capacity_estimate}
    \end{equation}
    Now we prove the dimension bound. For any map $V$ satisfying \eqref{lipschitz} and any sets $M\subset V N$,
    \begin{equation}
        n(M,q\varepsilon)\leq n(N,\varepsilon)\,\exp\left(\mathcal{C}_{\varrho_{t,s}}(N;\varepsilon,\delta\varepsilon)\right),
    \label{N}
    \end{equation}
    where $q=2\eta+\delta<1$. Choose a decreasing sequence $\{s_k\}_{k\ge0}$ with $s_0=S$ and $s_{k+1}\le s_k$. Then the estimates \eqref{iii} hold for each map $\mathscr{V}(s_k,s_{k+1})$, since all $s_k\le S$. By condition $(i)$,
    \begin{equation*}
        \mathcal{A}(s_k)\subset \mathscr{V}(s_k,s_{k+1})\mathcal{A}(s_{k+1}),\qquad 
        \mathcal{A}(t)\subset \mathscr{V}(t,S)\mathcal{A}(S).
    \end{equation*}
    Apply \eqref{N} to $V=\mathscr{V}(t,S)$ with $N=\mathcal{A}(S)$ and $M=\mathcal{A}(t)$:
    \begin{equation}
        n(\mathcal{A}(t),q\varepsilon)\leq n(\mathcal{A}(S),\varepsilon)\,\exp\left(\mathcal{C}_{\varrho_{t,S}}(\mathcal{A}(S);\varepsilon,\delta\varepsilon)\right).
    \label{7}
    \end{equation}
    For each $k\ge0$, apply \eqref{N} to $V=\mathscr{V}(s_k,s_{k+1})$ with $N=\mathcal{A}(s_{k+1})$ and $M=\mathcal{A}(s_k)$:
    \begin{equation}
        n(\mathcal{A}(s_k),q\varepsilon)\leq n(\mathcal{A}(s_{k+1}),\varepsilon)\,\exp\left(\mathcal{C}_{\varrho_{s_k,s_{k+1}}}(\mathcal{A}(s_{k+1});\varepsilon,\delta\varepsilon)\right).
    \label{8}
    \end{equation}
    Iterating \eqref{8} for $k=0,1,\dots,m-1$ and combining with \eqref{7} yields
    \begin{equation}
        n(\mathcal{A}(t),q^{m+1}\varepsilon)\leq n(\mathcal{A}(s_m),\varepsilon)\,\exp\left(\sum_{j=0}^{m}\mathcal{C}_j\right),
    \label{9}
    \end{equation}
    where
    \begin{equation*}
        \mathcal{C}_0=\mathcal{C}_{\varrho_{t,S}}(\mathcal{A}(S);\varepsilon,\delta\varepsilon),\qquad
        \mathcal{C}_{j+1}=\mathcal{C}_{\varrho_{s_j,s_{j+1}}}(\mathcal{A}(s_{j+1});q^{j}\varepsilon,\delta q^{j}\varepsilon).
    \end{equation*}
    By \eqref{capacity_estimate}, each $\mathcal{C}_j\leq \ln m_0\!\left(\frac{2K\sqrt{1+L^2}}{\delta}\right)$. Hence the sum is at most $(m+1)\ln m_0\!\left(\frac{2K\sqrt{1+L^2}}{\delta}\right)$.
    Since $q<1$, for any sufficiently small $\varepsilon>0$ we can choose $m$ such that $q^{m+1}\varepsilon$ is arbitrarily small; moreover $m+1\le \frac{\ln(1/\varepsilon)}{\ln(1/q)}+1$. Taking logarithms in \eqref{9} and letting $\varepsilon\to0$ gives
    \begin{equation*}
        \dim_f^X\mathcal{A}(t)\le \frac{\ln m_0\!\left(\frac{2K\sqrt{1+L^2}}{\delta}\right)}{\ln(1/q)}.
    \end{equation*}
    To have more explicit bounds, we can choose $\delta=\frac{1-2\eta}{2}$. Then $q=\frac{1+2\eta}{2}$ and $\frac{2}{\delta}=\frac{4}{1-2\eta}$. Substituting into above yields
    \begin{equation*}
        \dim_f^X\mathcal{A}(t)\leq \left[\ln\frac{2}{1+2\eta}\right]^{-1}\ln m_0\!\left(\frac{4K\sqrt{1+L^2}}{1-2\eta}\right)<\infty.
    \end{equation*}
    Since $\eta,K,L$ can be taken arbitrarily close to $\tilde\eta_t,\tilde K_t,\tilde L_t$, the bound is finite. Hence $\dim_f^{X}\mathcal{A}(t)<\infty$ for every fixed $t$. 
\end{proof}

For the rest of this section, we will prove all the assumptions in Theorem \ref{new_prop} for the attractors $\mathcal{A}_{\mu}(t)$ obtained in the previous sections  and obtain that they have finite fractional dimensions. We first impose an additional assumption on the external force
\begin{equation}\label{add_assump}
    \normm{g}:= \sup_{t\in \mathbb R}(\|g(t)\|_{L^2(\R)} + \|g'(t)\|_{L^2(\R)}) < +\infty.
\end{equation}
Fix $\mu\in (0,1]$ and let  $u_{1}(t)$ and $ u_{2}(t)$ be two orbits on the attractor $\{\mathcal{A}_{\mu }(t)\}$ for all $t \geq \tau$. We set $U(t) = u_1(t) - u_2(t)$. It follows that
\begin{equation}
	\begin{cases}
		U_{t} + U_{xxx} +( U\bar{u})_x - \mu U_{xx} + \lambda U = 0,& x\in \mathbb R, t>\tau, \\
		U(\tau, x)  = u_1(\tau,x)-u_2(\tau,x),& x\in\mathbb R,
	\end{cases}
	\label{ffd3}
\end{equation}
where $\bar{u}=\dfrac{u_{1}+u_{2}}{2}$.
\begin{lemma}\label{lm1}
	Let $U$ be the solution of \eqref{ffd3}. There exist $k_0\ge1$ and a constant $c_1>0$, depending only on $\lambda$, $\mu$ and $\normm{g}$, such that for all $k\ge k_0$ and $t>\tau$,
	\begin{equation}
		\|U(t)\|^{2}_{H^{1}} \leq e^{-\lambda (t-\tau)}\|U(\tau)\|^{2}_{H^{1}}\ + c_1\int_{\tau}^{t} \int_{|x|\leq k} U^{2}(s,x) dx ds.
		\label{ffd11}
	\end{equation}
 %    Furthermore, for all $t \geq \tau$,
 %    \begin{equation}
	% 	\|U(t)\|_{H^{1}} \leq \beta^{-1/2}e^{{c_1(t-\tau)}/2}\|U(\tau)\|_{H^{1}}.
	% 	\label{ffd18}
	% \end{equation}
\end{lemma} 
\begin{proof}
From estimates \eqref{f3} and \eqref{twostar} in Sections \ref{sec:2} and \ref{sec:KdV-Burgers}, and thanks to \eqref{add_assump} we see that
\begin{equation*}
    \sup_{t\in\R}\|\bar{u}(t)\|_{H^3} \leq \mathscr{C}_0,
\end{equation*}
where $\mathscr{C}_0$ depends on $\lambda, \mu$ and $\normm{g}$. Multiplying \eqref{ffd3} with $(U-\beta U_{xx})$, where $\beta=\min\{1,\mu/\mathscr C_0\}$, and integrating over $\R$, we find
\begin{equation}
	\begin{aligned}
	\frac 12\frac{d}{dt}\int_{\R}\bra{\beta U_x^2 +
    U^2}dx+ \mu\int_{\R}\bra{\beta U_{xx}^2 + U_x^2}dx + \lambda\int_{\R}\bra{\beta U_x^2 + U^2}dx\\
    = \int_{\R}(U\bar u)_x(-U + \beta U_{xx})dx.
	\end{aligned}
	\label{ffd5}   
\end{equation}
By integration by parts, we have
\begin{equation}\label{b1}
\begin{aligned}
    -\int_{\R}U(U\bar{u})_xdx &= -\frac 12\int_{\R}U^2\bar{u}_xdx \leq \frac 12\|\bar{u}_x\|_{L^{\infty}(|x|\geq k)}\int_{\R}U^2dx + \frac{1}{2}\mathscr{C}_0\int_{|x|\leq k}U^2dx\\
    &\leq \frac{\lambda}{4}\int_{\R}U^2dx + \frac{1}{2}\mathscr{C}_0\int_{|x|\leq k}U^2dx.
\end{aligned}
\end{equation}
provided $k$ is sufficiently large. For the second term on the right of \eqref{ffd5},
\begin{equation}\label{b2}
\begin{aligned}
    &\beta\int_{\R}U_{xx}(U\bar{u})_xdx = -\frac{\beta}{2}\int_{\R}U_x^2\bar{u}_xdx + \beta\int_{\R}U_{xx}U\bar{u}_xdx,\\
    &\leq \frac{\beta}{2} \mathscr{C}_0\int_{\R}U_{x}^2dx + \frac{\beta\mu}{2}\int_{\R}U_{xx}^2dx + \frac{\beta}{2\mu}\int_{\R}U^2|\bar{u}_x|^2dx,\\
    &\leq \frac{\beta}{2} \mathscr{C}_0\int_{\R}U_{x}^2dx + \frac{\beta\mu}{2}\int_{\R}U_{xx}^2dx + \frac{\beta  \mathscr{C}_0^2}{2\mu}\int_{|x|\leq k}U^2dx + \frac{\beta}{2\mu}\|\bar{u}_x\|_{L^{\infty}(|x|\geq k)}^2\int_{\R}U^2dx,\\
    &\leq \frac{\beta}{2}\mathscr{C}_0\int_{\R}U_{x}^2dx + \frac{\beta\mu}{2}\int_{\R}U_{xx}^2dx + \frac{\beta  \mathscr{C}_0^2}{2\mu}\int_{|x|\leq k}U^2dx + \frac{\lambda}{4}\int_{\R}U^2dx,
\end{aligned}
\end{equation}
for large enough $k$. Since $\beta\mathscr C_0\le\mu$, the terms $\frac\beta2\mathscr C_0\int U_x^2$ and $\frac{\beta\mu}2\int U_{xx}^2$ are absorbed by the second term on the left of \eqref{ffd5}.
By inserting \eqref{b2} and \eqref{b1} into \eqref{ffd5}, we obtain
\begin{equation}\label{b3}
    \frac{d}{dt}\int_{\R}(U^2 + \beta U_x^2)dx + \lambda \int_{\R}(U^2 + \beta U_x^2)dx \leq c_1\int_{|x|\leq k}U^2dx,
\end{equation}
where $c_1=\mathscr{C}_0 + \beta \mathscr{C}_0^2/\mu \leq 2\mathscr{C}_0$. By using the equivalent norm $\|f\|_{H^1}^2 \sim \|f\|_{L^2}^2 + \beta\|\partial_xf\|_{L^2}^2$, and applying Gronwall's lemma, we get
\begin{equation}\label{b4}
    \|U(t)\|_{H^1}^2 \leq e^{-\lambda (t-\tau)}\|U(\tau)\|_{H^1}^2 + c_1\int_{\tau}^t\int_{|x|\leq k}U^2(s,x)dxds.
\end{equation}
% Since $\beta\|U\|^2_{H^1} \leq \int (U^2+\beta U_x^2)dx \leq \|U\|^2_{H^1}$, \eqref{b4} yields the first estimate \eqref{ffd11}.\\
% \\To prove the second estimate, we return to the \eqref{b3}. Since
% \begin{equation*}
%     \int_{|x|\le k} U^2dx \leq \int_{\R} U^2dx \leq \int_{\R} (U^2+\beta U_x^2)dx,
% \end{equation*}
% and the damping term on the left of \eqref{b3} is nonnegative, we get
% \begin{equation*}
%     \frac{d}{dt}\int_{\R}\left(U^2 + \beta U_x^2\right)dx \leq c_1\int_{\R} (U^2+\beta U_x^2)dx.
% \end{equation*}
% Gronwall's lemma again gives
% \begin{equation*}
%     \int_{\R}\left(U^2(t) + \beta U_x^2(t)\right)dx \leq e^{c_1(t-\tau)}\int_{\R} \left(U^2(\tau)+\beta U_x^2(\tau)\right)dx,
% \end{equation*}
% and taking square roots, by the equivalence of norms,
% \begin{equation*}
%     \|U(t)\|_{H^{1}} \leq \beta^{-1/2}e^{c_1(t-\tau)/2}\|U(\tau)\|_{H^{1}}.
% \end{equation*}
\end{proof}

Let $k\geq1$. We consider the cutoff function $\varphi_{k} \in C^{\infty}(\R)$, which satisfies $0 \leq \varphi_k \leq 1,\quad \varphi_k = 1$ for $|x| \leq k, \varphi_k = 0$ for $|x| \geq 2k$, and
    \begin{equation}
    	|\partial_x^j\varphi_{k}(x)| \leq C_jk^{-j}, \quad j=0,1,2,3.
    \label{cutoff}
    \end{equation}
Denote by $U_{k}=\varphi_{k}U$. Then $U_k$ satisfies
\begin{equation}
    \left\{
	\begin{aligned}
		&(U_{k})_{t} + (U_{k})_{xxx} + (U_{k}\bar{u})_{x} + \lambda U_{k} -\mu (U_{k})_{xx} = m&,\\
		&U_{k}(\tau, x) = \varphi_{k}U(\tau,x)&,
	\end{aligned}
    \right.
	\label{ffd12}
\end{equation}
where
\begin{equation}\label{def_m}
	m = 3\partial_{x}(\partial_{x}U\partial_{x}\varphi_{k}) + U\partial_{x}^{3}\varphi_{k} + U\bar{u}\partial_{x}\varphi_{k} - 2\mu\partial_{x}U\partial_{x}\varphi_{k}-\mu U\partial_{x}^{2}\varphi_{k}.
\end{equation}
\begin{lemma}\label{lm2}
 There exist constants $c_6, c_7$ such that for all $t \geq \tau$,
\begin{equation}\label{f2}
    \int_{\tau}^t\|U_k(s)\|_{L^2}^2ds \le c_6\|U_k(\tau)\|_{L^2}^2 + \frac{c_7}{k^2}\|U(\tau)\|_{H^1}^2.
	%\|U_{k}(t)\|_{L^{2}} \leq e^{C(t-\tau)}\left(\|U_{k}(\tau)\|_{L^{2}} + k^{-1}\|U({\tau})\|_{H^1}\right).
\end{equation}
\end{lemma}

\begin{proof}
    Since $\bar{u}$ is bounded in $H^3(\R)$, we have $\|\bar{u}\|_{L^\infty} + \|\partial_x\bar{u}\|_{L^\infty} \leq 2\mathscr{C}_0$. 
    % Multiplying \eqref{ffd3} by $\varphi_k$, and using
    % \begin{align*}
    %     \varphi_kU_{xxx} &= (U_k)_{xxx} -3\partial_x(U_x\partial_x\varphi_k)-U\partial_x^3\varphi_k,&\\
    %     \varphi_k(U\bar{u})_x &= (U_k\bar{u})_x - U\bar{u}\partial_x\varphi_k,&\\
    %     \varphi_kU_{xx} &= (U_k)_{xx} - 2U_x\partial_x\varphi_k - U\partial_x^2\varphi_k,
    % \end{align*}
    % we obtain \eqref{ffd12}. 
    Multiplying \eqref{ffd12} by $U_{k}$ and integrating over $\R$, we get
	\begin{equation*}
		\frac{1}{2}\frac{d}{dt}\|U_k\|^2_{L^2} + \lambda \|U_k\|^2_{L^2} + \mu \|(U_k)_x\|^2_{L^2}= -\frac{1}{2} \int_{\R} \bar{u}_{x}U_k^{2} dx + \int_{\R} m U_k\ dx.
	\end{equation*}
    The first term on the right is bounded by $C\|U_k\|^2_{L^2}$. From \eqref{def_m}, we write $$\int_\R mU_kdx =  I_1 + I_2 + I_3 + I_4 + I_5$$ and estimate each term separately. Integrating by parts and using Young's inequality and \eqref{cutoff},
    \begin{equation*}
		\begin{aligned}
        I_1 = 3\int_\R \partial_x\left(\partial_x U\partial_x\varphi_k\right)U_k dx &= -3 \int_\R U_x\partial_x\varphi_k(U_k)_x dx, &\\
        &\leq \frac{\mu}{2} \|(U_k)_x\|^2_{L^2} + \frac{9}{2\mu}\int_{\R}|U_x|^2|\partial_x\varphi_k|^2 dx,&\\
        &\leq \frac{\mu}{2} \|(U_k)_x\|^2_{L^2} + \frac{C}{k^2}\|U\|^2_{H^1}.
        \end{aligned}
	\end{equation*}
    For the remaining terms, we use \eqref{cutoff} and Young's inequality:
    \begin{align*}
    I_2 &= \int_{\R} U \partial_x^3 \varphi_k U_k dx \le Ck^{-3} \|U\|_{L^2} \|U_k\|_{L^2} \le \frac{\lambda}{8} \|U_k\|_{L^2}^2 + C k^{-6} \|U\|_{L^2}^2, \\
    I_3 &= \int_{\R} U \bar{u} \partial_x \varphi_k U_k dx \le Ck^{-1} \|\overline{u}\|_{L^\infty} \|U\|_{L^2} \|U_k\|_{L^2} \le \frac{\lambda}{8} \|U_k\|_{L^2}^2 + C k^{-2} \|U\|_{L^2}^2, \\
    I_4 &= -2\mu \int_{\R} \partial_x U \partial_x \varphi_k U_k dx \le 2C\mu k^{-1} \|\partial_x U\|_{L^2} \|U_k\|_{L^2} \le \frac{\lambda}{8} \|U_k\|_{L^2}^2 + C k^{-2} \|U\|_{H^1}^2, \\
    I_5 &= -\mu \int_{\R} U \partial_x^2 \varphi_k U_k dx \le C\mu k^{-2} \|U\|_{L^2} \|U_k\|_{L^2} \le \frac{\lambda}{8} \|U_k\|_{L^2}^2 + C k^{-4} \|U\|_{L^2}^2.
    \end{align*}
    Since $k \geq1$, we have $k^{-4}, k^{-6} \le k^{-2}$, hence
    \begin{equation*}
        \int_{\R} m U_k dx \leq \frac{\mu}{2} \|(U_k)_x\|_{L^2}^2 + \frac{\lambda}{2} \|U_k\|_{L^2}^2 + \frac{C}{k^2} \|U\|_{H^1}^2.
    \end{equation*}
    Substituting this above, absorbing $\frac{\mu}{2} \|(U_k)_x\|_{L^2}^2$ into $\mu\|(U_k)_x\|^2_{L^2}$ on the left, we get
    \begin{equation*}
        \frac{d}{dt} \|U_k\|_{L^2}^2 + \lambda\|U_k\|_{L^2}^2 \le  \frac{C_4}{k^2} \|U(t)\|_{H^1}^2.
    \end{equation*}
    Integrating this on $(\tau,t)$ yields
    \begin{equation}\label{f1}
        \int_{\tau}^t\|U_k(s)\|_{L^2}^2 \le \lambda^{-1}\|U_k(\tau)\|_{L^2}^2 + \frac{C_5}{k^2}\int_{\tau}^t\|U(s)\|_{H^1}^2ds.
    \end{equation}
    At the same time, by integrating \eqref{b3} on $(\tau,t)$, we get, again using $\|f\|_{H^1}^2 \sim \|f\|_{L^2}^2 + \beta\|\partial_xf\|_{L^2}^2$,
    \begin{equation*}
        \lambda\int_{\tau}^t\|U(s)\|_{H^1}^2ds \le \|U(\tau)\|_{H^1}^2 + c_1\int_{\tau}^t\|U_k(s)\|_{L^2}^2ds
    \end{equation*}
    Inserting this one into \eqref{f1}, we get
    \begin{equation*}
        \int_{\tau}^t\|U_k(s)\|_{L^2}^2 \le \lambda^{-1}\|U_k(\tau)\|_{L^2}^2 + \frac{C_6}{k^2}\|U(\tau)\|_{H^1}^2 + \frac{C_7}{k^2}\int_{\tau}^t\|U_k(s)\|_{L^2}^2ds.
    \end{equation*}
    By choosing $k$ sufficiently large, we obtain the desired estimate \eqref{f2}.
    % \vskip 1in
    % By \eqref{ffd18}, $\|U(s)\|_{H^1}^2 \leq \beta^{-1}e^{c_1(s-\tau)} \|U(\tau)\|_{H^1}^2$ (\bao{maybe we should use \eqref{b4} here for a better control of the constant}). Gronwall's inequality on $[\tau, t]$ gives
    % \begin{equation*}
    %     \|U_k(t)\|_{L^2}^2 \leq e^{C_3(t-\tau)} \|U_k(\tau)\|_{L^2}^2 + \frac{C_4}{\beta k^2} \|U(\tau)\|_{H^1}^2 \int_\tau^t e^{C_3(t-s)} e^{c_1(s-\tau)} ds.
    % \end{equation*}
    % The integral is bounded by $(t-\tau)e^{Z(t-\tau)}$ with $Z = \max(C_3, c_1)$. Since $\frac{C_4}{\beta}(t-\tau) \leq e^{C_4(t-\tau)/\beta}$, we can absorb this factor into the exponent, and for a larger constant $\hat{C}$;
    % \begin{equation*}
    %     \|U_k(t)\|_{L^2}^2 \leq e^{\hat{C}(t-\tau)} \left(\|U_k(\tau)\|_{L^2}^2 + k^{-2} \|U(\tau)\|_{H^1}^2 \right).
    % \end{equation*}
    % Taking the square roots and using $\sqrt{a+b} \leq \sqrt{a} + \sqrt{b}$, we conclude
    % \begin{equation*}
    %     \|U_k(t)\|_{L^2} \le e^{C(t-\tau)} \left( \|U_k(\tau)\|_{L^2} + k^{-1} \|U(\tau)\|_{H^1} \right).
    % \end{equation*}
    \end{proof}
    \begin{theorem}
        For every $t \in \R$, the pullback attractor $\mathcal A_\mu(t)$ has a finite fractal dimension in $H^1(\R)$, i.e.,
        \begin{equation*}
            \dim_f^{H^1}\left(\mathcal A_\mu(t)\right) <+\infty.
        \end{equation*}
    \end{theorem}
    \begin{proof}
    We apply the Theorem \ref{new_prop} with
    \begin{equation*}
        X=H^1(\R),\qquad \mathcal A(t)=\mathcal{A}_\mu(t).
    \end{equation*}
    By Theorem \ref{compactness}, the family $\{\mathcal A(t)\}_{t\in\R}$ is compact in $H^1(\R)$. Let
    \begin{equation*}
        \mathscr V(t,s)=U_\mu(t,s)
    \end{equation*}
    be the process generated by \eqref{e1}. By the invariance of the pullback attractor, for $t\geq s$
    \begin{equation*}
        \mathcal A(t)=\mathscr V(t,s)\mathcal A(s),
    \end{equation*}
    so the condition $(i)$ of the Theorem \ref{new_prop} holds.\\
    \\Let $u_1,u_2\in\mathcal A(s)$ and set
    \begin{equation*}
        U(t)=\mathscr V(t,s)u_1-\mathscr V(t,s)u_2.
    \end{equation*}
    Then $U$ solves \eqref{ffd3}. From Lemma \ref{lm1}, for $k\geq k_0$,
    \begin{equation*}
        \|U(t)\|_{H^1}^2 \leq e^{-\lambda(t-s)}\|U(s)\|_{H^1}^2 + c_1\int_s^t\int_{|x|\le k}U^2(\tau,x)\,dx\,d\tau.
    \end{equation*}
    Since $\int_{|x|\leq k}U^2\,dx\leq \|U_k\|_{L^2}^2$, Lemma \ref{lm2} gives
    \begin{equation*}
        \int_s^t\|U_k(\tau)\|_{L^2}^2\,d\tau \leq c_6\|U_k(s)\|_{L^2}^2 + \frac{c_7}{k^2}\|U(s)\|_{H^1}^2.
    \end{equation*}
    Because $U_k(s)=\varphi_kU(s)$ and $\varphi_k=1$ on $|x|\leq k$, we have
    \begin{equation*}
        \|U_k(s)\|_{L^2}\leq \|U(s)\|_{L^2(|x|\leq 2k)}\leq \|U(s)\|_{H^1}.
    \end{equation*}
    Therefore,
    \begin{equation*}
    \|U(t)\|_{H^1}^2 \leq e^{-\lambda(t-s)}\|U(s)\|_{H^1}^2 + c_1c_6\|U(s)\|_{L^2(|x|\leq 2k)}^2 + \frac{c_1c_7}{k^2}\|U(s)\|_{H^1}^2.
    \end{equation*}
    Note that the constants $c_1, c_6, c_7$ are independent of $s, t,$ and $k$. Taking square roots and using $\sqrt{a+b}\le \sqrt a+\sqrt b$,
    \begin{equation*}
        \|U(t)\|_{H^1} \leq e^{-\lambda(t-s)/2}\|U(s)\|_{H^1} + (c_1c_6)^{1/2}\|U(s)\|_{L^2(|x|\leq 2k)} + \left(\frac{c_1c_7}{k^2}\right)^{1/2}\|U(s)\|_{H^1}.
    \end{equation*}
    Choose $k \geq k_0$ large enough so that
    \begin{equation*}
        \left(\frac{c_1c_7}{k^2}\right)^{1/2} \leq \frac{1}{8}.
    \end{equation*}
    Since $t$ is fixed, for $s$ sufficiently negative we also have
    \begin{equation*}
        e^{-\lambda(t-s)/2}\leq \frac{1}{8}.
    \end{equation*}
    Hence, for all such $s$,
    \begin{equation*}
        \|U(t)\|_{H^1} \leq \frac14\|U(s)\|_{H^1} + K\|U(s)\|_{L^2(|x|\leq 2k)},
    \end{equation*}
    where
    \begin{equation*}
        K=(c_1c_6)^{1/2}.
    \end{equation*}
    This is exactly the condition $(iii)$ of Theorem \ref{new_prop} with
    \begin{equation*}
        \eta_{t,s}=\frac14,\qquad n_1(v)=\|v\|_{L^2(|x|\le 2k)},\qquad n_2\equiv 0.
    \end{equation*}
    The semi-norm \(n_1\) is compact in \(H^1(\R)\), since the embedding $ H^1(|x|\le 2k)\hookrightarrow L^2(|x|\leq 2k)$ is compact. Moreover,
    \begin{equation*}
        \limsup_{s\to-\infty}\eta_{t,s}=\frac14<\frac12, \qquad \limsup_{s\to-\infty}K_{t,s}=K<+\infty.
    \end{equation*}
    For the condition $(ii)$, the same estimate gives
    \begin{equation*}
        \|U(t)\|_{H^1} \leq L_{t,s}\|U(s)\|_{H^1},
    \end{equation*}
    where
    \begin{equation*}
        L_{t,s}=\left(e^{-\lambda(t-s)}+c_1c_6+\frac{c_1c_7}{k^2}\right)^{1/2}.
    \end{equation*}
    Thus
    \begin{equation*}
        \limsup_{s\to-\infty}L_{t,s} = \left(c_1c_6+\frac{c_1c_7}{k^2}\right)^{1/2}<+\infty.
    \end{equation*}
    Therefore, all conditions of Theorem \ref{new_prop} are satisfied. We conclude that the pullback attractor $\mathcal{A}_\mu(t)$ has a finite fractal dimension in $H^1(\mathbb{R})$.
    \end{proof}
\vskip 0.5cm
\noindent{\bf Acknowledgments.} This work was carried out while the first
author was visiting the Mathematical Institute of the University of
G\"{o}ttingen. He would like to thank the Institute for its hospitality. N.A. and B.T. are supported by NAWI Graz.

\end{document}